\documentstyle[12pt]{article}
\newcommand{\sect}[1]{\section{#1}\setcounter{equation}{0}}

\font\mbn=msbm10 scaled \magstep1
\font\mbs=msbm7 scaled \magstep1
\font\mbss=msbm5 scaled \magstep1
\newfam\mbff
\textfont\mbff=\mbn
\scriptfont\mbff=\mbs
\scriptscriptfont\mbff=\mbss
\def\mbf{\fam\mbff}
\def\Re{{\mbf R}}

\def\Z{{\mbf Z}}
\def\Co{{\mbf C}}

\def\Di{{\mbf D}}
\newtheorem{Th}{Theorem}[section]
\newtheorem{Lm}[Th]{Lemma}
\newtheorem{C}[Th]{Corollary}

\newtheorem{Proposition}[Th]{Proposition}
\newtheorem{R}[Th]{Remark}
\newtheorem{E}[Th]{Example}
\author{Alexander Brudnyi\\
Department of Mathematics and Statistics\\
University of Calgary\\
Calgary, Canada}
\title{SOLVABLE MATRIX REPRESENTATIONS OF K\"{A}HLER GROUPS}
\date{}
\begin{document}
\maketitle
\begin{abstract}
{In the paper we study solvable matrix representations of fundamental
groups of compact K\"{a}hler manifolds ({\em K\"{a}hler groups}).
One of our main results is a factorization theorem for such representations.
As an application we prove that the universal covering of a compact 
K\"{a}hler manifold with a residually solvable fundamental group is 
holomorphically convex.}
\end{abstract}
\sect{\hspace*{-1em}. Introduction.} 
{\bf 1.1.} Let $M$ be a compact K\"{a}hler manifold. 
The fundamental group $\pi_{1}(M)$ of $M$ can be studied via its 
representations or equivalently in terms of locally constant sheaves on $M$. 
The set $Hom(\pi_{1}(M),GL_{n}(\Co))$ of all homomorphisms of 
$\pi_{1}(M)$ into the matrix Lie group $GL_{n}(\Co)$ has the natural 
structure of a complex affine algebraic variety. The group $GL_{n}(\Co)$ 
acts on $Hom(\pi_{1}(M),GL_{n}(\Co))$ by pointwise conjugation: 
$(gf)(s)=gf(s)g^{-1}$, $s\in\pi_{1}(M)$. The study of geometric
properties of $Hom(\pi_{1}(M),GL_{n}(\Co))$ is of interest because of the 
relation to the problem of classification of K\"{a}hler groups and 
description of complex manifolds which occur as nonramified coverings of a 
projective manifold. For any $\rho\in Hom(\pi_{1}(M),GL_{n}(\Co))$ let us 
denote by $X_{\rho}$ the union of all the irreducible components containing 
$\rho$. Let $D_{n}(\Co)\subset GL_{n}(\Co)$ be the subgroup of diagonal
matrices. The primary goal of this paper is to prove a factorization theorem 
for the elements of $X_{\rho}$ with $\rho\in Hom(\pi_{1}(M),D_{n}(\Co))$.
Before formulating the result let us first introduce some notation.

Let $N$ be a compact K\"{a}hler manifold. Hereafter
$\alpha_{N}:N\longrightarrow Alb(N)$ denotes the Albanese map of $N$ and
$N_{a}$ is the normalization of the image $\alpha_{N}(N)$. We also
denote by $\alpha_{N}^{n}:N\longrightarrow N_{a}$ the holomorphic map which
covers $\alpha_{N}$. 
\begin{Th}\label{te1}
Let $M$ be a compact K\"{a}hler manifold. There is a finite Galois 
covering $r:M_{H}\longrightarrow M$ with an abelian Galois group $H:=H(M)$ 
such that for any $\xi\in X_{\rho}$ with $\rho\in Hom(\pi_{1}(M),D_{n}(\Co))$
we have
$$
Ker(\alpha_{M_{H *}}^{n})\subset Ker(\xi)\ .
$$
\end{Th}
\begin{R}\label{r-1}
{\rm If $\rho\in Hom(\pi_{1}(M),D_{n}(\Co))$ is the trivial homomorphism,
then our method of proof gives also 
$Ker(\alpha_{M *}^{n})\subset Ker(\xi)$ for any $\xi\in X_{\rho}$.}
\end{R}
\begin{E}\label{ex1}
{\rm Let $T_{n}(\Co)\subset GL_{n}(\Co)$ be the subgroup of 
complex upper triangular matrices. For any 
$\rho\in Hom(\pi_{1}(M),T_{n}(\Co))$ let
$d(\rho)\in Hom(\pi_{1}(M),D_{n}(\Co))$ be the representation defined by the
diagonal of $\rho$. For the holomorphic diagonal matrix
$A(z):=diag[1,z,z^{2},...,z^{n-1}]$, $z\in\Co$,  we set 
$\rho(z):=A(z)^{-1}\rho A(z)$. Then $\rho(z)$ determines a holomorphic
deformation of $\rho$ to $\rho(0)=d(\rho)$. In particular, 
$\rho\in X_{d(\rho)}$.}
\end{E}
We now generalize the above example.
Let $V_{\rho}$ be a flat vector bundle on $M$ constructed by 
$\rho\in Hom(\pi_{1}(M),GL_{n}(\Co))$. Clearly, $V_{\rho}$
can be also considered as a holomorphic vector bundle. Assume that

(a)\ As a holomorphic bundle, $V_{\rho}$ is isomorphic to a
bundle $V_{\rho}^{\cal O}$ with structure group $T_{n}(\Co)$.\\
Let $Gr^{*}V_{\rho}^{\cal O}$ be the associated graded vector bundle with
cocycle defined as the diagonal of the cocycle of $V_{\rho}^{\cal O}$.
Then $Gr^{*}V_{\rho}^{\cal O}$ is the direct sum of rank-1 holomorphic
vector bundles on $M$.
Assume that

(b)\  The first Chern classes of components of the decomposition of
$Gr^{*}V_{\rho}^{\cal O}$ belong to $Tor(H^{2}(M),\Z)$.
\\
Denote the class of representations $\rho$ satisfying (a) and (b)
by $S_{n}^{\cal O}(M)$. Clearly, $Hom(\pi_{1}(M),T_{n}(\Co))\subset
S_{n}^{\cal O}(M)$. Let $D_{n}^{u}\subset D_{n}(\Co)$ be the subgroup of
unitary diagonal matrices and $p_{t}:M_{t}\longrightarrow M$ be
the Galois covering of $M$ with the Galois group 
$Tor(\pi_{1}(M)/D\pi_{1}(M))$. Here $DG=D^{1}G$ is the derived subgroup of
a group $G$. (We also set $D^{i}G=DD^{i-1}G$.) Under the
assumptions of Theorem \ref{te1} we prove
\begin{Th}\label{subset}
For any $\xi\in S_{n}^{\cal O}(M)$ there is
$\rho\in Hom(\pi_{1}(M),D_{n}^{u})$ such that
$p_{t}^{*}(\xi)\in X_{p_{t}^{*}(\rho)}$. Moreover, $r^{*}(\xi)$
factors through $\alpha_{M_{H}}$, that is,
$$
\xi\circ r_{*}=\xi'\circ\alpha_{M_{H *}}^{n}\ 
$$
for some $\xi'\in Hom(\pi_{1}(M_{H a}),GL_{n}(\Co))$.
\end{Th}
\begin{R}\label{r0}
{\rm Assume that $\xi\in S_{n}^{\cal O}(M)$ is such that
$Gr^{*}V_{\xi}^{\cal O}$ is the direct sum of $C^{\infty}$-trivial 
rank-1 holomorphic vector bundles (equivalently, all first Chern classes in
(b) above equal 0). Then we will show that $\xi\in X_{\rho}$ for some
$\rho\in Hom(\pi_{1}(M),D_{n}^{u})$.}
\end{R}
Let $E_{n}(M)$ be the class of
homomorphisms $\rho:\pi_{1}(M)\longrightarrow D_{n}(\Co)$ whose
diagonal characters $\rho_{ii}$ satisfy $\rho_{ii}=\exp(\widetilde\rho_{ii})$
for some $\widetilde\rho_{ii}\in Hom(\pi_{1}(M),\Co)$. By $S_{n}(M)$ we
denote the class of representations $\rho\in Hom(\pi_{1}(M),T_{n}(\Co))$ such
that the homomorphism $d(\rho)$ defined by the diagonal of $\rho$ 
belongs to $E_{n}(M)$. If $\pi:\widetilde{T_{n}(\Co)}\longrightarrow
T_{n}(\Co)$ is the universal covering, then, according to [O, Th.5.5], 
$\rho\in S_{n}(M)$ if and only if $\rho=\pi\circ\widetilde\rho$ for some 
$\widetilde\rho\in Hom(\pi_{1}(M),\widetilde{T_{n}(\Co)})$. In Lemma
\ref{coinc} we will show that $X_{\rho}\subset X_{d(\rho)}$ for any
$\rho\in S_{n}(M)$. Let $T_{n}^{u}\subset T_{n}(\Co)$ be the
subgroup of matrices with unitary diagonals and
$\widetilde{T_{n}^{u}}\subset\widetilde{T_{n}(\Co)}$ be the universal
covering of $T_{n}^{u}$.
We set
$S_{n}^{u}(M):=S_{n}(M)\cap Hom(\pi_{1}(M),T_{n}^{u})$ and define
$$
G_{s}^{0}(M)=\bigcap_{n\geq 1}\bigcap_{\rho\in S_{n}(M)}Ker(\rho)\ \ \ \ \
G_{su}^{0}(M)=\bigcap_{n\geq 1}\bigcap_{\rho\in S_{n}^{u}(M)}Ker(\rho).
$$
Notice that the above groups are kernels of ${\cal P}$-completion 
homomorphisms $\eta^{\cal P}:\pi_{1}(M)\longrightarrow\pi_{1}(M)^{\cal P}$
where ${\cal P}$ is the category of groups $\widetilde{T_{n}(\Co)}$ and
$\widetilde{T_{n}^{u}}$, respectively.
In general we can say only that $G_{s}^{0}(M)\subseteq G_{su}^{0}(M)$. 
However, if $M$ is a compact K\"{a}hler manifold then 
\begin{Th}\label{eq}
$G_{s}^{0}(M)=G_{su}^{0}(M)$.
\end{Th}
{\bf 1.2.} We apply Theorem \ref{te1} to the problem of the
holomorphic convexity of a Galois covering of a compact K\"{a}hler manifold.
Central to the subject is a conjecture of Shafarevich according to which the 
universal covering $M_{u}$ of a projective manifold $M$ is holomorphically 
convex, meaning that for every infinite sequence of points without limit 
points in $M_{u}$ there exists a holomorphic function unbounded on this 
sequence. In recent years many new powerful methods have been developed to 
investigate the structure of the fundamental groups and universal coverings 
of complex projective manifolds. These methods lead to many new remarkable 
results. In particular quite a few positive results on the Shafarevich's 
conjecture were proved, e.g. the results of F. Campana, H. Grauert, 
R. Gurjar, L. Katzarkov, M. Nori, M. Ramachandran, C. Simpson, S. Shasrty, 
K. Zuo, S.T. Yau 
(see e.g. [C], [Ka1,2], [KaR], [Ko1,2], [N], [S], [Y]). It was also shown by 
F.Bogomolov and L.Katzarkov [BK] that a counterexample to the
Shafarevich conjecture exists provided that one can answer affirmatively a
certain group theoretic question.

For a compact K\"{a}hler manifold $M$ we define
\begin{equation}\label{resid}
G_{s}(M):=\bigcap_{n\geq 1}\bigcap_{\rho\in Hom(\pi_{1}(M),T_{n}(\Co))}
Ker(\rho)
\end{equation}
and let $M_{s}$ be the regular covering of $M$ with the transformation
group $\pi_{1}(M)/G_{s}(M)$.
\begin{Th}\label{te2}
$M_{s}$ is a holomorphically convex manifold. Thus if
$p:M_{1}\longrightarrow M$ is a finite unbranched covering of $M$ such
that the family $\{Hom(\pi_{1}(M_{1}),T_{n}(\Co))\}_{n\geq 1}$
separates the elements of $\pi_{1}(M_{1})$ then the universal covering 
$M_{u}$ of $M$ is holomorphically convex.
\end{Th}
As a corollary we obtain the result of Katzarkov [Ka1] 
on the holomorphic convexity of the Malcev covering of any smooth projective 
manifold.
\begin{C}\label{stein}
Let $M$ be a compact K\"{a}hler manifold such that $G_{s}(M)=\{e\}$ and
$\pi_{2}(M)=0$. Then $M_{u}$ is a Stein manifold.
\end{C}
{\bf 1.3.} We formulate several corollaries on the structure of 
K\"{a}hler groups.

Let $C$ be a smooth compact complex curve of genus $g\geq 1$. Any character 
$\xi\in Hom(\pi_{1}(C),\Co^{*})$ can be written as $\exp(\widetilde\xi)$ for 
some $\widetilde\xi\in Hom(\pi_{1}(C),\Co)$. In turn, $\widetilde\xi$ is 
defined by integration along the paths of a harmonic
1-form $\omega$ on $C$. Let $\omega=\omega'+\omega''$ be the type 
decomposition of $\omega$ into the sum of holomorphic and antiholomorphic
1-forms and $\widetilde\xi''\in Hom(\pi_{1}(C),\Co)$ be the homomorphism 
defined by integration of $\omega''-\overline{\omega''}$. We set
$\xi_{u}''=\exp(\widetilde\xi'')$. Then $\xi_{u}''\in Hom(\pi_{1}(C),U_{1})$
where $U_{1}$ is the unitary group. Notice also that $\xi_{u}''$ is
uniquely defined by $\xi$.
Let $f:M\longrightarrow C$ be a holomorphic surjective map with connected
fibres of a compact K\"{a}hler manifold $M$ onto a smooth compact complex
curve $C$ of genus $g\geq 1$.
\begin{Th}\label{transf}
Assume that $\rho\in Hom(\pi_{1}(M),T_{n}(\Co))$ satisfies\\
$$
\begin{array}{l}
\displaystyle
(1)\ \ d(\rho)=\bigoplus_{i=1}^{n}f^{*}(\rho_{i})\ \ {\rm for\ some\ }\ \
\rho_{i}\in Hom(\pi_{1}(C),\Co^{*});\\
\\
\displaystyle
(2)\ \ \rho_{iu}''\otimes(\rho_{ju}'')^{-1}\ \ {\rm is\ not\ a\ torsion\ 
character\ for\ any}\ \ i\neq j.
\end{array}
$$
Then there is $\widetilde\rho\in Hom(\pi_{1}(C),T_{n}(\Co))$ such that
$\rho=f^{*}(\widetilde\rho)$.
\end{Th}
\begin{R}\label{antiho}
{\rm If $\rho\in Hom(\pi_{1}(M),T_{n}^{u})$ then (2) acquires
the form $\rho_{i}\otimes\rho_{j}^{-1}$ is not torsion for any
$i\neq j$.}
\end{R}
Let $G$ be the fundamental group of a compact Riemann surface of genus
$g\geq 2$ and $G_{1}$ be an extension of $G$ by $\Z^{k}$, $k\geq 1$,
such that the conjugate action $s$ of $G$ on $\Z^{k}$ is trivial on $DG$. 
Thinking of $\Z^{k}$ as a subgroup of $\Z^{k}\otimes\Co$ we can extend
$s$ to a representation $s':G\longrightarrow GL_{k}(\Co)$ so that
$s'(g)|_{\Z^{k}}=s(g)$ for any $g\in G$. Then $DG\subset Ker(s')$ and so
$s'$ admits a decomposition $\oplus_{j=1}^{m}s_{j}$ where $s_{j}$ is
equivalent to a nilpotent representation $G\longrightarrow T_{k_{j}}(\Co)$
with diagonal character $\rho_{j}$.
\begin{Th}\label{group}
Assume that a K\"{a}hler group $F$ is defined as
$$
\{e\}\longrightarrow K\stackrel{i}{\longrightarrow} F
\stackrel{j}{\longrightarrow} G_{1}\longrightarrow
\{e\}\ 
$$
and $j^{-1}(\Z^{k})\subset F$ does not admit a surjective 
homomorphism onto a free group with countable number of generators. Then all 
characters $\rho_{j}$ are torsion.
\end{Th}
The above theorem produces many examples of non-K\"{a}hler groups even in
the case when $K$ is finitely presented. For instance, a semidirect product
of $G$ with $\Z^{k}$ defined by $s:G\longrightarrow SL_{k}(\Z)$ as above
such that not all $\rho_{j}$ are torsion is not a K\"{a}hler group. 

We will complete this section by
\begin{Th}\label{kob}
Let $M$ be a compact K\"{a}hler manifold satisfying the conditions of
Corollary \ref{stein} and $M_{H}$ be the same as in Theorem \ref{te1}.
Assume also that
$$
dim_{\Co}M\geq\frac{1}{2}rank(\pi_{1}(M)/D^{2}\pi_{1}(M))\ .
$$
Then
$dim_{\Co}M=\frac{1}{2}rank(\pi_{1}(M)/D^{2}\pi_{1}(M))$ and
$\pi_{1}(M_{H})$ is a free abelian group.\\
If, in addition,
$$
\pi_{i}(M_{u})=0\ {\rm for}\  
1\leq i\leq dim_{\Co}M
$$ 
then $\alpha_{M_{H}}:M_{H}\longrightarrow Alb(M_{H})$ is biholomorphic.
\end{Th}
\begin{R}\label{rtfn}
{\rm If $\pi_{1}(M)$ is residually torsion free nilpotent satisfying
$\pi_{2}(M)=0$ and $dim_{\Co}M\geq \frac{1}{2}rank(\pi_{1}(M)/D\pi_{1}(M))$ 
then $dim_{\Co}M=\frac{1}{2}rank(\pi_{1}(M)/D\pi_{1}(M))$ and 
$\pi_{1}(M)$ is free abelian. If, in addition, $\pi_{i}(M_{u})=0$ for 
$1\leq i\leq dim_{\Co}M$ then $\alpha_{M}$ is biholomorphic.}
\end{R}
In particular, any compact complex manifold $M$ with a nonzero
Kobayashi semi-metric satisfying the conditions of 
Theorem \ref{kob} or Remark \ref{rtfn} is not K\"{a}hler.

In the next section we reduce Theorem \ref{te1} to the
case of representations from $Hom(\pi_{1}(M),T_{n}(\Co))$.  Also we
describe a classification theorem for a class of flat connections on
$M\times\Co^{n}$ with triangular (0,1)-components (see [Br]) which is
crucial in the proof of Theorem \ref{te1}.
\sect{\hspace*{-1em}. Flat Connections and
Flat Vector Bundles.}
{\bf 2.1.}
In this section we collect some results on the relation between 
flat vector bundles and the flat connections which determine them. 

Let $M$ be a compact K\"{a}hler manifold and $gl_{n}(\Co)$ be the Lie
algebra of $GL_{n}(\Co)$. Consider the family of $C^{\infty}$-trivial flat 
vector bundles on $M$ of complex rank $n$. It is well-known that every bundle
from this family is determined by a flat connection on the trivial bundle 
$M\times\Co^{n}$, that is, by a $gl_{n}(\Co)$-valued differential
1-form $\omega$ on $M$ satisfying
\begin{equation}\label{eq1}
d\omega-\omega\wedge\omega=0\ .
\end{equation}
Denote the class of flat connections by ${\cal A}_{n}(M)$. The group
$C^{\infty}(M,GL_{n}(\Co))$ acts by d-gauge transforms on the set
${\cal A}_{n}(M)$:
\begin{equation}\label{eq2}
d_{g}(\alpha)=g^{-1}\alpha g-g^{-1}dg,\ \ \ g\in C^{\infty}(M,GL_{n}(\Co)),
\ \ \alpha\in {\cal A}_{n}(M)\ .
\end{equation}
Denote the corresponding quotient set by ${\cal B}_{n}(M)$. (We regard
${\cal B}_{n}(M)$ as the set of d-gauge equivalent classes of connections
from ${\cal A}_{n}(M)$.) Let $p:M_{u}\longrightarrow M$ be the universal
covering of $M$ and $p^{*}(\omega)$ be the pullback of $\omega$. Fix a point 
$x\in M_{u}$. Then there is a unique solution 
$I_{x}(\omega)\in C^{\infty}(M_{u},GL_{n}(\Co))$ of the equation
$dF=p^{*}(\omega)F$ satisfying $I_{x}(\omega)(x)=E_{n}$ (here $E_{n}$ denotes
the unit matrix). $I_{x}(\omega)$ can be obtained by the Picard method of 
successive approximations (also called  {\em iterated path integration}).
Let $|\cdot|$ denote the $l_{2}$-norm in the space of $n\times n$ matrices.
Then we also have 
\begin{equation}\label{est1}
|I_{x}(\omega)(y)|\leq e^{cr(y)|\omega|_{M}}
\end{equation}
where $r$ is the distance from $x$ in the metric pulled back from $M$, 
$c:=c(M)>0$ is a
constant and $|\cdot|_{M}$ is an $L_{2}$-norm on the space
of $gl_{n}(\Co)$-valued differential 1-forms on $M$. Moreover, if the form
$\omega=\omega(z)$ depends holomorphically on a parameter $z$ varying in the 
unit disk $\Di\subset\Co$, then $I_{x}(\omega(z))(y)$ is holomorphic in $z$
for any $y\in M_{u}$.

Let $o(x):=\{h(x)\}_{h\in\pi_{1}(M)}$ be the orbit of 
$x$ under the action of the
covering group $\pi_{1}(M)$. Then the restriction $I_{x}(\omega)|_{o(x)}$ 
determines a homomorphism of $\pi_{1}(M)$ into $GL_{n}(\Co)$:
$$
I_{x}(\omega)((h_{1}h_{2})(x))=I_{x}(\omega)(h_{1}(x))\cdot
I_{x}(\omega)(h_{2}(x)),\ \ \ h_{1},h_{2}\in\pi_{1}(M)\ .
$$
If $x,y\in M_{u}$ are distinct points then there is a matrix 
$C(x,y)\in GL_{n}(\Co)$ such that
$I_{y}(\omega)|_{o(y)}=C(x,y)^{-1}I_{x}(\omega)|_{o(x)}C(x,y)$.
Further, for any $g\in C^{\infty}(M,GL_{n}(\Co))$
there is a matrix $C\in GL_{n}(\Co)$ such that 
$$
I_{x}(d_{g}(\alpha))=C^{-1}I_{x}(\alpha)C, \ \ \ \alpha\in {\cal A}_{n}(M)\ .
$$ 
Thus
$I_{x}$ determines an injective map ${\cal I}_{x}:
{\cal B}_{n}(M)\rightarrow Hom(\pi_{1}(M),GL_{n}(\Co))/GL_{n}(\Co)$. Its 
image consists of equivalence classes of homomorphisms whose associated flat 
vector bundles are $C^{\infty}$-trivial. 

Another way to construct a representation by a flat connection 
$\omega\in {\cal A}_{n}(M)$ is as follows 
(see also, e.g. [O, Sec.5,6]).
Let $(U_{i})_{i\in I}$ be an open covering of $M$ by
sets diffeomorphic to Euclidean balls and let 
$f_{i}\in C^{\infty}(M,GL_{n}(\Co))$ be a solution of $df=\omega f$ on
$U_{i}$. If we set $c_{ij}=f_{i}^{-1}f_{j}$ on $U_{i}\cap U_{j}$,
then $\{c_{ij}\}$ is a
locally constant cocycle and so it determines a flat vector bundle $E(\omega)$
on $M$. Then a standard procedure (see, e.g. [KN, Ch.2, Sec.9])
allows to construct a $\rho\in Hom(\pi_{1}(M),GL_{n}(\Co))$ by $E(\omega)$.
Notice that $\rho$ coincides with $I_{x}(\omega)|_{o(x)}$ for a 
suitable $x\in M_{u}$.

Let $E$ be a $C^{\infty}$-trivial flat vector bundle on $M$.
Consider the flat vector bundle $End(E)$. Let $\{c_{ij}\}$ be a locally
constant cocycle on an open acyclic covering $(U_{i})_{i\in I}$ of $M$
determining $E$.
A family $\{\eta_{i}\}_{i\in I}$ of matrix-valued $p$-forms satisfying
\begin{equation}\label{vf}
\eta_{j}=c_{ij}^{-1}\eta_{i}c_{ij}\ \ \ \ {\rm on}\ \ \ U_{i}\cap U_{j}
\end{equation}
is, by definition, a $p$-form with values in the bundle $End(E)$.
According to (\ref{vf}) the operators $d$ and $\wedge$ are well-defined
on the set of matrix-valued $p$-forms. In particular, one can consider
$gl_{n}(\Co)$-valued 1-forms $\alpha$ with values in $End(E)$ 
satisfying the flatness condition $d\alpha-\alpha\wedge\alpha=0$.
Let $h$ be a linear  $C^{\infty}$-automorphism of $E$ determined by 
the family $\{h_{i}\}_{i\in I}$, $(h_{i}\in C^{\infty}(U_{i},GL_{n}(\Co)))$, 
satisfying
$$
h_{j}=c_{ij}^{-1}h_{i}c_{ij}\ \ \ \ {\rm on}\ \ \ U_{i}\cap U_{j}\ .
$$
Then a $d$-gauge transform $d_{h}^{E}$ defined on the set of 
$gl_{n}(\Co)$-valued 1-forms $\alpha$ with values in $End(E)$ is 
given by the formula
$$
d_{h}^{E}(\alpha)=h^{-1}\alpha h-h^{-1}dh\ .
$$
Clearly, $d_{h}^{E}$ preserves the class of 1-forms satisfying
the flatness condition. Let $\psi$ be a flat connection on $M\times\Co^{n}$
determining $E$ and $g_{i}\in C^{\infty}(U_{i},GL_{n}(\Co))$ be a
solution of $df=\psi f$ such that $g_{i}^{-1}g_{j}=c_{ij}$. Consider
the map $\tau_{\psi}$ defined on the set of matrix-valued 1-forms
$\alpha$ on $M$ by
\begin{equation}\label{map}
\tau_{\psi}(\alpha_{i})=g_{i}^{-1}(\alpha_{i}-\psi_{i})g_{i},
\ \ \ \alpha_{i}:=\alpha|_{U_{i}},\ \ \ \psi_{i}:=\psi|_{U_{i}}\ .
\end{equation}
Then $\tau_{\psi}$ maps the set ${\cal A}_{n}(M)$ isomorphically onto
the set of $End(E)$-valued matrix 1-forms satisfying the flatness
condition.
Further, if $h\in C^{\infty}(M,GL_{n}(\Co))$ then the family 
$\{g_{i}^{-1}hg_{i}\}_{i\in I}$ determines an element of $Aut(E)$.
In what follows we identify $C^{\infty}(M,GL_{n}(\Co))$ with $Aut(E)$.
Then we have
$$
\tau_{\psi}\circ d_{h}=d_{h}^{E}\circ\tau_{\psi}\ \ \
{\rm for\ every}\ \ \ h\in C^{\infty}(M,GL_{n}(\Co))\ 
$$
(for the proof, see e.g. [Br, Prop.2.2]). 
\\
{\bf 2.2.}  In this part we describe some results of [Br] (see also
[BO] for the general case).
Let $\omega$ be a flat connection on $M\times\Co^{n}$
such that the $(0,1)$-component $\omega''$ of $\omega$ is an upper
triangular matrix form. Doing a $d$-gauge transform with a diagonal
matrix-function we may assume without loss of generality that
$diag(\omega'')$ is a matrix-valued harmonic $(0,1)$- form on $M$, in 
particular, that it is $d$-closed ($M$ is K\"{a}hler). Consider another flat 
connection on $M\times\Co^{n}$ defined by the system of ODEs
\begin{equation}\label{diago}
df=\psi f\
\end{equation}
where $\psi:=diag(\omega'')-diag(\overline{\omega''})$.
Let $\{f_{i}\}_{i\in I}$, $f_{i}\in C^{\infty}(M,D_{n}^{u})$, be a family of 
local solutions of (\ref{diago}) defined on an open acyclic covering 
$(U_{i})_{i\in I}$ of $M$, where
$D_{n}^{u}\subset U_{n}$ is the Lie subgroup of unitary diagonal $n\times n$ 
matrices. Then 
the cocycle 
$$ 
c_{ij}:=f_{i}^{-1}f_{j}\ \ \ {\rm on}\ \ \  U_{i}\cap U_{j} 
$$
determines a flat vector bundle $V_{\psi}$ which is the direct sum of 
$C^{\infty}$-trivial complex rank-1 flat vector bundles with unitary structure
group. In what follows $V_{\psi}$-valued harmonic forms are defined with 
respect to Laplacians constructed by the natural flat Hermitian metric on 
$End(V_{\psi})$. Consider now 
the map $\tau_{\psi}$ defined in Section 2.1 with $\psi$ as above.
Clearly $\tau_{\psi}(\omega)$ is an $End(V_{\psi})$-valued 1-form  with
a nilpotent (0,1)-component. Recall that a family $\{\eta_{i}\}_{i\in I}$
of matrix-valued 1-forms satisfying
$$
\eta_{j}=c_{ij}^{-1}\eta_{i}c_{ij}\ \ \ {\rm on}\ \ \ U_{i}\cap U_{j}
$$
is a nilpotent 1-form with values in $End(V_{\psi})$ if every $\eta_{i}$
takes its values in the Lie algebra of the Lie group of upper triangular
unipotent matrices. It was proved in [Br, Prop.2.3] that

{\bf (A)} {\em There is a $g\in C^{\infty}(M,T_{n}(\Co))$
such that $\tau_{\psi}(d_{g}(\omega))=\eta_{1}+\eta_{2}$,
where the (1,0)-form $\eta_{1}$ is $\partial$-closed and $\eta_{2}$ is a
$d$-closed nilpotent antiholomorphic 1-form. In particular, if the form
$\omega$ is upper triangular then $\eta_{1}$ is also upper triangular.}\\
This implies that $\eta_{1}$ can be decomposed into the sum
$\alpha+\partial h$, where $\alpha$ is its harmonic component in the Hodge
decomposition. Then the Hodge decomposition and the 
$\partial\overline\partial$-lemma (see e.g. [ABCKT, Lm.7.39]) imply

{\bf (B)} {\em $\alpha$ is $d$-harmonic and $(\alpha+\eta_{2})\wedge
(\alpha+\eta_{2})$ represents 0  in the de Rham cohomology group
$H^{2}(M, End(V_{\psi}))$.}

Conversely, let $\alpha$ be an $End(V_{\psi})$-valued $d$-harmonic 
$(1,0)$-form and let $\theta$ be a $d$-harmonic 
nilpotent $End(V_{\psi})$-valued $(0,1)$-form. Then we have
(see [Br, Prop.2.4]):

{\bf (C)} {\em Let $(\alpha+\theta)\wedge(\alpha+\theta)$ represent 0 in
$H^{2}(M,End(V_{\psi}))$. Then there exists a unique (up to a flat
additive summand) section $h$ such that $\alpha+\partial h+\theta$
satisfies the flatness condition. In particular, there is a flat connection
$\omega$ on $M\times\Co^{n}$ with an upper triangular $(0,1)$-component
such that $\tau_{\psi}(\omega)=\alpha+\partial h+\theta$. If, in 
addition, $\alpha$ is upper triangular then $\omega$ is also upper triangular.
}\\
The following examples will be used in the sequel.
\begin{E}\label{ex2}
{\rm 1. Let $C$ be a compact complex curve of genus $g\geq 1$. 
Any representation
$\rho\in Hom(\pi_{1}(C),T_{n}(\Co))$ determines a $C^{\infty}$-trivial
flat vector bundle on $C$. Let $d(\rho)\in Hom(\pi_{1}(C),D_{n}(\Co))$
be the representation defined by the diagonal of $\rho$. Then $d(\rho)$
is defined by the system of ODEs $df=\phi f$ with diagonal harmonic 1-form
$\phi$. Let $\phi=\phi'+\phi''$ be the type decomposition. Denote
by $d(\rho)_{u}''\in Hom(\pi_{1}(C),D_{n}^{u})$ the representation defined
by equation $df=(\phi''-\overline{\phi''})f$. Let $V_{d(\rho)_{u}''}$
be the flat vector bundle associated to $d(\rho)_{u}''$. Then 
$V_{d(\rho)_{u}''}=\oplus_{i=1}^{n}V_{i}$ where
each $V_{i}$ is a complex rank-1 flat vector bundle on $C$ with 
structure group $U_{1}$. Further, $End(V_{d(\rho)_{u}''})=\oplus_{i,j=1}^{n}
V_{i}^{*}\otimes V_{j}$. Assume that $V_{i}^{*}\otimes V_{j}$ is not
trivial for any $i\neq j$. Then it is to see (e.g. by duality) that the
de Rham cohomology group $H^{2}(C,V_{i}^{*}\otimes V_{j})$ is zero. According
to the above results the equivalence class of $\rho$ is defined by
an upper triangular harmonic 1-form $\eta$ with values in 
$End(V_{d(\rho)_{u}''})$ such that $\eta\wedge\eta$ represents 0 in
$H^{2}(C,End(V_{d(\rho)_{u}''}))$ and the (0,1)-component of $\eta$ is 
nilpotent. Notice that the identity $H^{2}(C,V_{i}^{*}\otimes V_{j})=0$ for 
$i\neq j$ 
implies that the first condition for $\eta$ is always fulfilled. Conversely,
any upper-triangular $End(V_{d(\rho)_{u}''})$-valued harmonic 1-form 
$\eta$ with a nilpotent (0,1)-component determines the equivalence class of a 
representation from $Hom(\pi_{1}(C),T_{n}(\Co))$.\\
2. Let $M$ be a compact K\"{a}hler manifold, $\rho\in S_{n}^{\cal O}(M)$ and 
$p_{t}:M_{t}\longrightarrow M$ be the Galois covering with the
Galois group $Tor(\pi_{1}(M)/D\pi_{1}(M)$. Let $V_{\rho}$ be the
associated flat vector bundle. Then according to condition (b) in the 
definition of $S_{n}^{\cal O}(M)$, the first Chern classes of components of
the decomposition of $Gr^{*}(p_{t}^{*}V_{\rho}^{\cal O})$ are zero.
Therefore $p_{t}^{*}V_{\rho}$ is defined by a flat connection on
$M_{t}\times\Co^{n}$ with a triangular (0,1)-component and the results of
this section describe the equivalence class of $p_{t}^{*}V_{\rho}$.}
\end{E}
{\bf 2.3.} In this section we reduce Theorem \ref{te1} to the
case of solvable matrix representations.

Let $p_{1}(g_{1},...,g_{k})=e,...,p_{l}(g_{1},...,g_{k})=e$ be a finite 
presentation of $\pi_{1}(M)$. Here $g_{1},...,g_{k}\in\pi_{1}(M)$ is a set 
of generators. Then $Hom(\pi_{1}(M),GL_{n}(\Co))$ is defined 
(up to polynomial isomorphisms) as an affine algebraic subvariety
of $(GL_{n}(\Co))^{k}$ by equations
$$
p_{1}(X_{1},...,X_{k})=E_{n},...,p_{l}(X_{1},...,X_{k})=E_{n},\ \ \ \ 
X_{1},...,X_{k}\in GL_{n}(\Co)\ .
$$
Let $\rho\in Hom(\pi_{1}(M),GL_{n}(\Co))$. Assume that there is an 
irreducible component $X\subset Hom(\pi_{1}(M),GL_{n}(\Co))$ of $\rho$ 
such that $dim_{\Co}X\geq 1$.
Since $X$ is a complex algebraic variety of pure dimension $\geq 1$, 
resolution of singularities shows that there is a non-constant
holomorphic map $f:\Di\longrightarrow X$ 
such that $f(0)=\rho$. Let ${\cal F}_{\rho}$ be the family of all such maps.
Then there is an open neighbourhood $O\subset X$ of
$\rho$ such that for any $w\in O$ there is an $f\in {\cal F}_{\rho}$ so that
$w\in f(\Di)$. Any $f\in {\cal F}_{\rho}$ determines a holomorphic
deformation of $\rho$ (which we denote by the same letter).
In particular, we have
$f(z)(g)=\sum_{i=0}^{\infty}f_{i}(g)z^{i}$, 
$g\in\pi_{1}(M)$, such that $f(z)\in Hom(\pi_{1}(M),GL_{n}(\Co))$,
$f_{0}=\rho$, $f_{i}$ are functions on 
$\pi_{1}(M)$ with values in $gl_{n}(\Co)$ and the series converges uniformly 
on each compact subset of $\Di$ for any $g\in\pi_{1}(M)$.

Let $\Co[[z]]$ be the ring of formal power series over $\Co$. Let
$I_{m}\subset\Co[[z]]$ denote the ideal generated by $z^{m}$. Let
$\pi_{m}:\Co[[z]]\longrightarrow\Co[[z]]/I_{m}=:A_{m}[z]$ be the quotient 
homomorphism. For any $h\in\Co[[z]]$ we identify $\pi_{m}(h)$ with
the first $m$ terms of the expansion of $h$. Let
$GL_{n}(A_{m}[z])$ be the group of invertible
$n\times n$ matrices with entries from $A_{m}[z]$. Any $a\in GL_{n}(A_{m}[z])$
can be written as $a=\sum_{i=0}^{m-1}a_{i}z^{i}$, where 
$a_{i}\in gl_{n}(\Co)$. The element $a$ is invertible if and only if
$a_{0}\in GL_{n}(\Co)$. Similarly, denote by $gl_{n}(A_{m}[z])$ the algebra
of all $n\times n$ matrices with entries from $A_{m}[z]$. 
The homomorphism $\pi_{m}$ induces a similar homomorphism (denoted
by the same letter) $GL_{n}(\Co)[[z]]\longrightarrow GL_{n}(A_{m}[z])$ of 
the space of formal matrix power series 
which associate to each series the first $m$ terms of its Taylor expansion.
Assume now that $f\in {\cal F}_{\rho}$ is a homomorphism depending 
holomorphically on $z\in\Di$. Then 
$\pi_{m}(f)\in Hom(\pi_{1}(M),GL_{n}(A_{m}[z]))$.
Obviously we have
\begin{Proposition}\label{kern}
$$
\bigcap_{z\in\Di}Ker(f(z))=
\bigcap_{i\geq 0}Ker(\pi_{i}(f(z)))\ \ \ \ \ \ \  \Box
$$
\end{Proposition}
The main result of this section is the following reduction in Theorem 
\ref{te1}. We retain the notation introduced in Theorem \ref{te1}.
\begin{Proposition}\label{triang}
Theorem \ref{te1} is a consequence of the following result:\\
For any $\rho\in Hom(\pi_{1}(M),T_{n}(\Co))$, $n\geq 1$,
we have $Ker(\alpha_{M_{H *}}^{n})\subset Ker(\rho)$. 
\end{Proposition}
{\bf Proof.}
Let $S_{n}^{d}\subset gl_{n}(\Co)[[z]]$ be the subalgebra of matrix power 
series $a=\sum_{i=0}^{\infty}a_{i}z^{i}$ with diagonal $a_{0}$. We
set $F_{m,n}:=\pi_{m}(S_{n}^{d})\subset gl_{n}(A_{m}[z])$. Any
element $b\in F_{m,n}$ admits an expansion $b=\sum_{i=0}^{m-1}b_{i}z^{i}$
with diagonal $b_{0}$. Further, $F_{m,n}$ is a finite-dimensional Lie 
algebra over $\Co$ with the bracket: $[a,b]:=ab-ba$, $a,b\in F_{m,n}$. 
Let $F_{m,n}^{*}:=\{b\in F_{m,n}:\ 
b=\sum_{i=0}^{m-1}b_{0}z^{i}, \ b_{0}\in D_{n}\}$ denote the 
group of invertible elements of $F_{m,n}$ (with respect to multiplication). 
Then $F_{m,n}^{*}$ is a complex solvable Lie group. Let $\phi:
F_{m,n}^{*}\longrightarrow GL_{p}(\Co)$, $p=dim_{\Co}F_{m,n}$,
be the faithful representation obtained by the left 
multiplication of the elements of $F_{m,n}$ by the elements of $F_{m,n}^{*}$.
We set $V_{k}:=\{v\in F_{m,n} :\ v=\sum_{i=k}^{m-1}v_{i}z^{i}\}$. 
Then
$$
F_{m,n}=V_{0}\supset V_{1}\supset V_{2}\supset...\supset V_{m}=\{0\}
$$
is a chain of invariant $F_{m,n}^{*}$-submodules with respect to $\phi$.  
The quotient representation $F_{m,n}^{*}\longrightarrow GL(V_{i+1}/V_{i})$ is
defined by left multiplication of matrices from $gl_{n}(\Co)$ by 
invertible diagonal matrices. In particular, $\phi$ is equivalent to a 
representation $\widetilde\phi:F_{m,n}^{*}\longrightarrow T_{p}(\Co)$ 
(defined by the choice of a suitable basis in $F_{m,n}$).

Consider now $\rho\in Hom(\pi_{1}(M),D_{n}(\Co))$. Let $f\in {\cal F}_{\rho}$
be a holomorphic deformation of $\rho$. Then $\pi_{m}(f)\in
Hom(\pi_{1}(M),F_{m,n}^{*})$. Further, consider 
$f_{m}:=\widetilde\phi\circ\pi_{m}(f(z))$. It is easy to see that the 
diagonal characters of $f_{m}$ are the same as for $\rho$. 
Clearly $f_{m}\in Hom(\pi_{1}(M),T_{p}(\Co))$. Assume now that
Theorem \ref{subset} is valid for any 
$f_{m}$, $m\geq 0$, and any $n\geq 1$. Then
according to Proposition \ref{kern}, $Ker(f(z))$ contains
$Ker(\alpha_{M_{H *}}^{n})$ for any $z\in\Di$. 

Let $X_{1}\subset X_{\rho}$ be an irreducible component of
$Hom(\pi_{1}(M),GL_{n}(\Co))$ containing $\rho$ and $O\subset X_{1}$
be an open neighbourhood of $\rho$ such that for any
$w\in O$, $w\neq\rho$, there is $f\in {\cal F}_{\rho}$
satisfying $w\in f(\Di)$. Then 
the assumption of the proposition and the previous arguments applied to 
such $f(z)$ imply that
$Ker(\alpha_{M_{H *}}^{n})\subset Ker(w)$ for any $w\in O$. Let
$g\in Ker(\alpha_{M_{H *}}^{n})$ and $g=p(g_{1},...,g_{k})$ be a word 
representing $g$. Here $g_{1},...,g_{k}\in\pi_{1}(M)$ are generators.
Then the complex algebraic matrix function $p(X_{1},...,X_{k})$, $X_{1},...,
X_{k}\in GL_{n}(\Co)$, equals 0 on an open set $O\subset X_{1}$. Since
$X_{1}$ is irreducible, $p\equiv 0$ on $X_{1}$. This
argument applied to each $g\in Ker(\alpha_{M_{H *}}^{n})$ shows that
$Ker(\alpha_{M_{H *}}^{n})\subset Ker(q)$ for any $q\in X_{1}$. The
same is valid for other irreducible components of $X_{\rho}$ and for
any $\rho\in Hom(\pi_{1}(M),D_{n}(\Co))$.
This implies the required statement.
\ \ \ \ \ $\Box$
\sect{\hspace*{-1em}. Fibering K\"{a}hler Manifolds.}
In this section we study some properties of fundamental groups of compact
K\"{a}hler manifolds which admit surjective holomorphic maps
with connected fibers onto compact Riemann surfaces. We also prove
Theorems \ref{transf} and \ref{group}.
\\
{\bf 3.1.} Let $f:M\longrightarrow C$ be a holomorphic surjective map with
connected fibres of a compact K\"{a}hler manifold $M$ onto a smooth compact
complex curve $C$. Consider a flat vector bundle $L$ on $C$ of complex rank 
1 with unitary structure group and let $E=f^{*}L$ be the pullback of this 
bundle to $M$. Let $\omega\in\Omega^{1}(E)$
be a holomorphic 1-form with values in $E$ 
(by the Hodge decomposition $\omega$ is $d$-closed).
\begin{Lm}\label{vanish}
Assume that $\omega|_{V_{z}}=0$ for the fibre $V_{z}:=f^{-1}(z)$ over a 
regular value $z\in C$. Then $\omega|_{V}=0$ for any fibre $V$ of $f$.
\end{Lm}
\begin{R}
{\rm If $V$ is not smooth the lemma asserts that $\omega$ equals 0 on the 
smooth part of $V$.}
\end{R}
{\bf Proof.} Denote by $S$ 
the set of non-regular values of $f$. Then $S\subset C$ consists of a finite 
number of points. First consider fibres of $f$ over regular values.
By Sard's theorem, $f:M\setminus f^{-1}(S)\longrightarrow C\setminus S$ is
a fibre bundle with connected fibres. 
According to our assumption there is a 
fibre $V_{z}$ of this bundle such that $\omega|_{V_{z}}=0$. 
Then $\omega|_{V}=0$ for any fibre $V$ of $f|_{M\setminus f^{-1}(S)}$.
In fact, for any fibre $V$ there is an open neighbourhood $O_{V}$ of $V$
diffeomorphic to $\Re^{2}\times V$ such
that $E|_{O_{V}}$ is the trivial flat vector bundle (here we have 
used the fact
that $E|_{V}$ is trivial because it is pullback of a bundle defined over a 
point). Any $d$-closed holomorphic 1-form defined on $O_{V}$ that vanishes on
$V$ is $d$-exact, and so it vanishes on each fibre contained in $O_{V}$. 
Starting with a tubular neighbourhood $O_{V_{z}}$ and taking into account 
that $\omega|_{O_{V_{z}}}$ can be considered as a $d$-closed holomorphic 
1-form (because $E|_{O_{V_{z}}}$ is trivial) we obtain
that $\omega|_{V}=0$ for any fibre $V\subset O_{V_{z}}$. Then the required
statement follows by induction if we cover $M\setminus f^{-1}(S)$ by open 
tubular neighbourhoods of fibres of $f$ and use the fact that 
$C\setminus S$ is connected.

Consider now fibres over the singular part $S$.  Let $O_{x}\subset C$ be the
neighbourhood of a point $x\in S$ such that 
$O_{x}\setminus S=O_{x}\setminus\{x\}$ is biholomorphic to 
$\Di\setminus\{0\}$. In particular, $\pi_{1}(O_{x}\setminus S)=\Z$. 
Without loss of generality we may assume that $V=f^{-1}(x)$ is a deformation 
retract of an open neighbourhood $U_{x}$ of $V$ and 
$W_{x}:=f^{-1}(O_{x})\subset U_{x}$ (such $U_{x}$ exists, e.g.,
by the triangulation theorem of the pair $(M,V)$, see [L, Th.2]). Since 
$E|_{V}$ is trivial, the bundle $E|_{U_{x}}$ is also trivial. Therefore we 
can regard $\omega$ as a $d$-closed holomorphic 1-form on $U_{x}$. Moreover, 
we already have proved that $\omega|_{F}=0$ for any fibre 
$F\subset W_{x}\setminus f^{-1}(S)$. So $\omega$ equals 0 on each fibre of 
the fibration $f:W_{x}\setminus f^{-1}(S)\longrightarrow O_{x}\setminus S$. 
This implies that there is a $d$-closed
holomorphic  1-form $\omega_{1}$ on $O_{x}\setminus S$ such that $\omega=
f^{*}(\omega_{1})$. Assume now that $\omega|_{V}\neq 0$. Then integration
of $\omega|_{W_{x}}$ along paths determines a non-trivial homomorphism 
$h_{\omega,x}:\pi_{1}(W_{x})\longrightarrow\Co$ whose image
is isomorphic to $\Z$. In fact, 
since embedding $W_{x}\setminus f^{-1}(S)\hookrightarrow W_{x}$ induces a 
surjective homomorphism of fundamental groups, we can integrate $\omega$ by 
paths contained in $W_{x}\setminus f^{-1}(S)$ obtaining the same image in 
$\Co$. Further, for any path $\gamma\subset W_{x}\setminus f^{-1}(S)$ we have
$\int_{\gamma}\omega=\int_{\gamma_{1}}\omega_{1}$, where $\gamma_{1}$ is a 
path in $O_{x}\setminus S$ representing the element 
$f_{*}(\gamma)\in\pi_{1}(O_{x}\setminus S)$. But 
$\pi_{1}(O_{x}\setminus S)=\Z$ and therefore 
$h_{\omega,x}(\pi_{1}(W_{x}))\cong\Z$. Recall that any path in $W_{x}$ is 
homotopically equivalent inside $U_{x}$ to a path contained in $V$ and so 
$\omega|_{V}$ determines a homomorphism $h$ of 
$\Gamma:=(\pi_{1}(V)/D\pi_{1}(V))/_{torsion}$ into $\Co$ with image 
isomorphic to $\Z$. Without loss of generality we may assume
that $V$ is smooth (for otherwise, we apply the arguments below to each
irreducible component of a desingularization of $V$). Integrating holomorphic
1-forms on $V$ along paths we embed  $\Gamma$ into some $\Co^{k}$ as a 
lattice of rank $2k$. Then there is a linear holomorphic functional 
$f_{\omega}$ on $\Co^{k}$ such that $h=f_{\omega}|_{\Gamma}$. Since 
$rk(h(\Gamma))=1$, $f_{\omega}$ equals 0 on a subgroup $H\subset\Gamma$ 
isomorphic to $\Z^{2k-1}$. In particular,
$f_{\omega}=0$ on the vector space $R=span(H)$ of real dimension $2k-1$. But
$Ker(f_{\omega})$ is a complex vector space. So $f_{\omega}=0$ on $\Co^{k}$ 
which implies that $\omega|_{V}=0$. This contradicts our assumption and 
proves the required statement for fibres over the points of $S$.

The lemma is proved.\ \ \ \ \ $\Box$
\\
{\bf 3.2.} Let $f: M\longrightarrow C$ be a surjective holomorphic map with
connected fibres of a compact K\"{a}hler manifold $M$ onto a smooth compact
complex curve $C$. Then the induced homomorphism $f_{*}:\pi_{1}(M)
\longrightarrow\pi_{1}(C)$ is surjective.
Denote by $T_{2}\subset T_{2}(\Co)$ the Lie group of matrices of the form
\[
\left(
\begin{array}{cc}
a&b\\ 0&1\\
\end{array}
\right)\ \ \ \ \ (a\in\Co^{*},\ b\in\Co) .
\]
For a homomorphism $\rho\in Hom(\pi_{1}(M),T_{2})$ we let $\rho_{a}\in
Hom(\pi_{1}(M),\Co^{*})$ denote the upper diagonal character of $\rho$.
\begin{Proposition}\label{normsubgr}
Let $\rho\in Hom(\pi_{1}(M),T_{2})$ be such that $\rho_{a}$ is a non-unitary
character. Assume that $\rho_{a}|_{Ker(f_{*})}$ is trivial. Then
$\rho|_{Ker(f_{*})}$ is also trivial.
\end{Proposition}
{\bf Proof.}
Let $p_{t}:M_{t}\longrightarrow M$ be the Galois covering with  
transformation group $Tor(\pi_{1}(M)/D\pi_{1}(M))$. Then 
$\rho|_{\pi_{1}(M_{t})}$ determines a $C^{\infty}$-trivial complex rank-2 
flat vector bundle on $M_{t}$, because 
$\rho_{a}|_{\pi_{1}(M_{t})}=\exp(\rho'|_{\pi_{1}(M_{t})})$ for some
$\rho'\in Hom(\pi_{1}(M),\Co)$. In particular, $\rho|_{\pi_{1}(M_{t})}$
can be defined by a flat connection 
\[
\Omega=\left(
\begin{array}{cc}
\omega&\eta\\0&0\\
\end{array}
\right);\ \ \ \ \ \ \  d\Omega-\Omega\wedge\Omega=0
\]
on $M_{t}\times\Co^{2}$ where $\omega$ is a $d$-harmonic 1-form on $M_{t}$
lifted from $M$.
Moreover, since $\rho_{a}|_{Ker(f_{*})}$ is trivial, $\omega$ is the
pullback by $f\circ p_{t}$ of a $d$-harmonic 1-form defined on $C$. Let
$\omega=\omega_{1}+\omega_{2}$ be the type decomposition of $\omega$ into the
sum of holomorphic and antiholomorphic 1-forms lifted from $C$. Denote by 
$E_{\rho}$ the rank-1 flat vector bundle on $M_{t}$ with unitary structure 
group constructed by the flat connection $\omega_{2}-\overline{\omega_{2}}$ 
and by $E_{0}=M_{t}\times\Co$ the trivial flat vector bundle. Observe that 
$E_{\rho}$ is pullback of a flat bundle on $C$. In particular, 
$E_{\rho}|_{V}$ is trivial for any fibre $V$ of $f\circ p_{t}$.  According 
to the results of Sections 2.1 and 2.2 the equivalence
class of $\rho|_{\pi_{1}(M_{t})}$ is determined by a $d$-harmonic 1-form 
$\theta$
with values in $End(E_{\rho}\oplus E_{0})$ satisfying $\theta\wedge\theta$
represents 0 in $H^{2}(M_{t},End(E_{\rho}\oplus E_{0}))$. More precisely,
\[
\theta=
\left(
\begin{array}{cc}
\omega_{1}&\eta'\\0&0\\
\end{array}
\right)
\]
where $\eta'$ is a $d$-harmonic 1-form with values in $E_{\rho}$ satisfying
$\omega_{1}\wedge\eta'$ represents 0 in $H^{2}(M_{t},E_{\rho})$.
Let $M_{t}\stackrel{g_{1}}{\longrightarrow}M_{1}\stackrel{g_{2}}
{\longrightarrow}C$ be the Stein factorization of 
$f\circ p_{t}$. Here $g_{1}$ is a morphism with connected
fibres onto a smooth curve $M_{1}$ and $g_{2}$ is a finite morphism. 
Then for any fibre $V\hookrightarrow
M_{t}$ of $g_{1}$ we have that $p_{t}|_{V}:V\longrightarrow p_{t}(V)$ is a
regular covering of the fibre $p_{t}(V)$ of $f$ with a finite abelian
transformation group.
\begin{Lm}\label{van}
$\eta'|_{V}=0$ for fibres over regular values of $g_{1}$.
\end{Lm}
Based on this statement we, first, finish the proof of the proposition and 
then will prove the lemma. From the lemma and Lemma \ref{vanish} it follows 
that $\eta'|_{W}=0$ for any fibre $i:W\hookrightarrow M_{t}$ of $g_{1}$.  Thus
$\rho|_{\pi_{1}(M_{t})}$ is trivial on $i_{*}(\pi_{1}(W))\subset
\pi_{1}(M_{t})$. But $\pi_{1}(W)$ is a subgroup of a finite index in
$\pi_{1}(p_{t}(W))$ and the image of $\rho|_{j_{*}(\pi_{1}(p_{t}(W)))}$
consists of unipotent matrices by the assumption of the proposition. (Here
$j:p_{t}(W)\hookrightarrow M$ is embedding.) Therefore 
$\rho|_{j_{*}(\pi_{1}(p_{t}(W)))}$ is trivial for any $W$.
Denote by $E$ the flat vector bundle on $M$ associated to $\rho$. Then 
we have proved that $E|_{W}$ is trivial for any fibre $W$ of $f$. Further, 
let $(U_{i})_{i\in I}$ be an open
covering of $C$ such that $W_{i}:=f^{-1}(U_{i})$ is an open neighbourhood of a
fibre $V_{i}=f^{-1}(x_{i})$, $x_{i}\in U_{i}$, and $W_{i}$ is deformable onto
$V_{i}$. Since $E|_{V_{i}}$ is trivial, $E|_{W_{i}}$ is also trivial. In
particular, $E$ is defined by a locally constant cocycle $\{c_{ij}\}$ 
defined on the covering $(W_{i})_{i\in I}$. But then $\{c_{ij}\}$ is the 
pullback of a cocycle defined on $(U_{i})_{i\in I}$ because the fibres of
$f$ are connected. This cocycle determines a bundle $E'$ on $C$ such that 
$f^{*}E'=E$. Then $\rho$ is the pullback of a homomorphism
$\rho'\in Hom(\pi_{1}(C),T_{2})$ constructed by $E'$.

This completes the proof of the proposition modulo Lemma \ref{van}.\\ 
{\bf Proof of Lemma \ref{van}.} 
Observe that $\omega_{1}\neq 0$ in the above definition of $\theta$ 
because $\rho_{a}$ is non-unitary. We also regard $\eta'\neq 0$, 
for otherwise, the image of $\rho$ consists of diagonal matrices and the 
required statement is trivial. Let
$i:V\hookrightarrow M_{t}$ be a fibre over a regular value of $g_{1}$. For
any $\lambda\in\Co^{*}$ consider the form
\[
\theta_{\lambda}=\left(
\begin{array}{cc}
\lambda\omega_{1}&\eta'\\0&0\\
\end{array}
\right).
\]
Clearly, $\theta_{\lambda}\wedge\theta_{\lambda}$ represents 0 in
$H^{2}(M_{t},End(E_{\rho}\oplus E_{0}))$ and thus, according to the results of
Section 2.2, it determines a
representation $\rho_{\lambda}\in Hom(\pi_{1}(M_{t}),T_{2})$ with the upper
diagonal character $\rho_{\lambda a}$ defined by the flat connection
$\lambda\omega_{1}+\omega_{2}$. In particular, the family
$\{\rho_{\lambda a}\}$ contains infinitely many different characters. 
Assume that $\rho|_{i_{*}(\pi_{1}(V))}$ is not trivial. Then it can be 
determined by the restriction  $\theta|_{i(V)}$. But according to our 
assumption
$\omega_{1}|_{i(V)}=0$ and $E_{\rho}|_{i(V)}$ is a trivial flat vector 
bundle. Thus $\rho|_{i_{*}(\pi_{1}(V))}$ is
defined by $\eta'|_{i(V)}\in H^{1}(i(V),\Co)$. Let $\psi\in
Hom(i(V),\Co)$ be a homomorphism obtained by integration of
$\eta'|_{i(V)}$ along paths generating $\pi_{1}(i(V))$. Then
\[
\rho|_{i_{*}(\pi_{1}(V))}=\left(
\begin{array}{cc}
1&\psi\\0&1\\
\end{array}
\right).
\]
We also obtain that $\rho_{\lambda}|_{i_{*}(\pi_{1}(V))}=\rho|_{
i_{*}(\pi_{1}(V))}$ and so it is not trivial for any $\lambda\in\Co^{*}$. 
We will prove now that the family $\{\rho_{\lambda a}\}$ contains finitely 
many different characters. This contradicts our assumption and shows that
$\rho|_{i_{*}(\pi_{1}(V))}$ is trivial and so $\eta'|_{V}=0$.

Let $H:=i_{*}(\pi_{1}(V))$ be a normal subgroup of $\pi_{1}(M_{t})$ (the 
normality follows from the exact homotopy sequence for the bundle defined
over regular values of $g_{1}$).
Denote by $H_{ab}$ the quotient $H/DH$. Then $H_{ab}$ is an 
abelian group of a finite rank. Moreover, $H_{ab}$ is a normal subgroup of 
$K:=\pi_{1}(M_{t})/DH$ and the group $K_{1}:=K/H_{ab}$ is finitely generated. 
Note that $\rho_{\lambda}$ induces a homomorphism 
$\widehat\rho_{\lambda}:K\longrightarrow T_{2}$ for any
$\lambda\in\Co^{*}$ with the same image as for $\rho_{\lambda}$
because the image of $\rho_{\lambda}|_{H}$ is abelian. In particular,
$\widehat\rho_{\lambda}(H_{ab})$ is a normal subgroup of $\widehat\rho(K)$. 
Since by our assumption $\widehat\rho_{\lambda}(H_{ab})$ is a non-trivial 
subgroup of unipotent matrices, from the identity
\[\left(
\begin{array}{cc}
c&b\\0&1\\
\end{array}
\right)\cdot
\left(
\begin{array}{cc}
1&v\\0&1\\
\end{array}
\right)\cdot
\left(
\begin{array}{cc}
c&b\\0&1\\
\end{array}
\right)^{-1}=
\left(
\begin{array}{cc}
1&c\cdot v\\0&1\\
\end{array}
\right)\ \ \ \ \ (c\in\Co^{*},\ b,\ v\in\Co)
\]
it follows that the action of $\widehat\rho_{\lambda}(K)$ on
$\widehat\rho_{\lambda}(H_{ab})$ by
conjugation is defined by multiplication of non-diagonal elements of
$\widehat\rho_{\lambda}(H_{ab})$ by elements of 
$\rho_{\lambda a}(\pi_{1}(M_{t}))$. Note also that $Tor(H_{ab})$ belongs to
$Ker(\widehat\rho_{\lambda})$ for any $\lambda$.
Thus if $H'\cong\Z^{s}$ is a maximal free 
abelian subgroup of $H_{ab}$ then $\rho_{\lambda a}(\pi_{1}(M_{t}))$ 
consists of 
eigen values of matrices from $SL_{s}(\Z)$ obtained by the natural action 
(by conjugation) of $K_{1}$ on $H'$. Since $K_{1}$ is
finitely generated, the number of different $\rho_{\lambda a}$ is finite
which is false.

This contradiction completes the proof of the lemma.\ \ \  \ \ $\Box$
\\
{\bf 3.3.} Let $f:M\longrightarrow C$ be a holomorphic map with connected 
fibres of a compact K\"{a}hler manifold $M$ onto a smooth
compact complex curve $C$ of genus $g\geq 1$. Then $\pi_{1}(M)$ is defined 
by the exact sequence
$$
\{e\}\longrightarrow Ker(f_{*})\longrightarrow\pi_{1}(M)\longrightarrow
\pi_{1}(C)\longrightarrow\{e\}\ . $$
Moreover, $Ker(f_{*})$ is a finitely generated group. Let 
$G\subset Ker(f_{*})$ be a normal subgroup of $\pi_{1}(M)$. Then the 
quotient group $R:=\pi_{1}(M)/G$ is
defined by the sequence 
$$ 
\{e\}\longrightarrow Ker(f_{*})/G\longrightarrow
R\longrightarrow\pi_{1}(C)\longrightarrow\{e\}\ . 
$$ 
Assume that

(1)\ \ $H:=Ker(f_{*})/G$ is a free abelian group of rank $k\geq 1$;

(2)\ \ the natural action of $D\pi_{1}(C)$ on $H$ is trivial.
\\
From (2) it follows that the action $s$ of $\pi_{1}(C)$ on $H$ determines an 
action of $\Z^{2g}\cong\pi_{1}(C)/D\pi_{1}(C)$ on $H$. Identifying $H$
with $\Z^{k}$ we can think of $H$ as a subgroup (lattice) of 
$Z^{k}\otimes\Co=\Co^{k}$. Then $s$ determines a representation 
$s':\pi_{1}(C)\longrightarrow GL_{k}(\Co)$ such that 
$Ker(s')\subset D\pi_{1}(C)$ and $s'(g)|_{H}=s(g)$ for any $g\in\pi_{1}(C)$. 
Since $s'$ descends to a representation $\Z^{2g}\longrightarrow GL_{k}(\Co)$,
it admits a decomposition $s'=\oplus_{j=1}^{m}s_{j}$ where $s_{j}$ is
equivalent to a nilpotent representation 
$\pi_{1}(C)\longrightarrow T_{k_{i}}(\Co)$ with a 
diagonal character $\rho_{j}$. Here $\sum_{j=1}^{m}k_{j}=k$.
\begin{Proposition}\label{curve}
All characters $\rho_{j}$ are torsion.
\end{Proposition}
{\bf Proof.} By definition, $R$ is an extension of $\pi_{1}(C)$ by $H$. It 
is well known, see, e.g. [G, Ch.I, Sec.6], that the class of extensions 
equivalent to $R$ is uniquely defined by an element 
$c\in H^{2}(\pi_{1}(C),H)$ where the cohomology is defined by the action
$s$ of $\pi_{1}(C)$ on $H$. Let $f\in Z^{2}(\pi_{1}(C),H)$ be a cocycle 
determining $c$. Then one can define a representative of the 
equivalence class of extensions as the direct product $H\times\pi_{1}(C)$ 
with 
multiplication 
$$
\begin{array}{c}
(h_{1},g_{1})\cdot (h_{2},g_{2})=
(h_{1}+s(g_{1})(h_{2})+f(g_{1},g_{2}),g_{1}\cdot g_{2});\\ 
\\
h_{1},h_{2}\in H,\ g_{1},g_{2}\in\pi_{1}(C)\ .
\end{array}
$$ 
The natural embedding $H\hookrightarrow\Z^{k}\otimes\Co (=\Co^{k})$
determines an embedding $i$ of $R$ into the group $R'$ defined as
$\Co^{k}\times\pi_{1}(C)$ with multiplication 
$$
\begin{array}{r}
(v_{1},g_{1})\cdot (v_{2},g_{2})=
(v_{1}+s'(g_{1})(v_{2})+f(g_{1},g_{2}),g_{1}\cdot g_{2});\\
\\
f\in Z^{2}(\pi_{1}(C),H),\ v_{1},v_{2}\in\Co^{k},\ g_{1},g_{2}\in\pi_{1}(C)\ .
\end{array}
$$ 
Here we regard $f$ as an element of $Z^{2}(\pi_{1}(C),\Co^{k})$ defined by
the action $s'$. From the decomposition $s'=\oplus_{j=1}^{m}s_{j}$ it follows
that there is an invariant $\pi_{1}(C)$-submodule $V_{j}\subset\Co^{k}$ of
$dim_{\Co}V_{j}=k-1$ such that $W_{j}=\Co^{k}/V_{j}$ is a one-dimensional
$\pi_{1}(C)$-module and the action of $\pi_{1}(C)$ on $W_{j}$ is defined as
multiplication by the character $\rho_{j}$. Then, by definition, $V_{j}$ is 
a normal subgroup of $R'$ and the the quotient group $R_{j}=R'/V_{j}$ is 
defined by the sequence 
$$
\{e\}\longrightarrow\Co\longrightarrow
R_{j}\longrightarrow\pi_{1}(C)\longrightarrow\{e\}\  .
$$ 
Here the action of $\pi_{1}(C)$ on $\Co$ is multiplication by the character 
$\rho_{j}$. Further, the equivalence class of extensions isomorphic to 
$R_{j}$ is defined by an element $c_{j}\in H^{2}(\pi_{1}(C),\Co)$ (the 
cohomology is defined by the above action of $\pi_{1}(\Co)$ on $\Co$). We
will assume that the character $\rho_{j}$ is non-trivial (for otherwise,
$\rho_{j}$ is clearly torsion). Let us
denote by $t_{j}$ the composite homomorphism $\pi_{1}(M)\longrightarrow
R\stackrel{i}{\longrightarrow}R'\longrightarrow R_{j}$.

Let $E_{\rho_{j}}$ be a complex rank-1 flat vector bundle on $C$ constructed 
by $\rho_{j}\in Hom(\pi_{1}(C),\Co^{*})$. Since $C$ is a 
$K(\pi_{1}(C),1)$-space, there is a natural isomorphism of the above group 
$H^{2}(\pi_{1}(C),\Co)$ and the \v{C}ech cohomology group 
$H^{2}(C, {\bf E_{\rho_{j}}})$ of the sheaf of locally
constant sections of $E_{\rho_{j}}$, see e.g. [M, Ch.1, Complement to Sec.2].
But each flat vector bundle on $C$ is $C^{\infty}$-trivial and each 
homomorphism from $Hom(\pi_{1}(C),\Co^{*})$ can be continuously deformed 
inside of $Hom(\pi_{1}(C),\Co^{*})$ to the trivial homomorphism. Therefore 
by the index theorem 
$$ 
dim_{\Co}H^{0}(C,{\bf E_{\rho_{j}}})-dim_{\Co}H^{1}(C,{\bf E_{\rho_{j}}})+
dim_{\Co} H^{2}(C,{\bf E_{\rho_{j}}})=\chi(\pi_{1}(C))=2-2g\ .
$$
Note that $H^{0}(C,{\bf E_{\rho_{j}}})=0$ because $\rho_{j}$ is non-trivial.
Furthermore, $H^{1}(C,{\bf E_{\rho_{j}}})$ is in a one-to-one correspondence 
with the set of non-equivalent representations 
$\rho:\pi_{1}(C)\longrightarrow T_{2}\subset T_{2}(\Co)$ with the
upper diagonal character $\rho_{j}$. Using the identity
$\prod_{i=1}^{g}[e_{i},e_{g+i}]=e$ for generators 
$e_{1},..,e_{2g}\in\pi_{1}(C)$ we easily obtain that 
$dim_{\Co}H^{1}(C,{\bf E_{\rho_{j}}})=2g-2$. Thus we have
$H^{2}(C,{\bf E_{\rho_{j}}})=0$. This shows that $H^{2}(\pi_{1}(C),\Co)=0$ and
$R_{j}$ is isomorphic to the semidirect product of $\Co$ and $\pi_{1}(C)$, 
i.e., $R_{j}=\Co\times\pi_{1}(C)$ with multiplication 
$$
(v_{1},g_{1})\cdot (v_{2},g_{2})=
(v_{1}+\rho_{j}(g_{1})\cdot v_{2},g_{1}\cdot g_{2}),\ \ \
v_{1},v_{2}\in \Co,\ g_{1},g_{2}\in\pi_{1}(C)\ .
$$ 
Let us determine a map $\phi_{j}$ of $R_{j}$ to $T_{2}$
by the formula
\[
\phi_{j}(v,g)=\left(
\begin{array}{cc}
\rho_{j}(g)&v\\ 0&1\\
\end{array}
\right)
\]
Obviously, $\phi_{j}$ is a correctly defined homomorphism with  upper 
diagonal character $\rho_{j}$. Hence 
$\phi_{j}\circ t_{j}:\pi_{1}(M)\longrightarrow T_{2}(\Co)$ is a homomorphism 
which is non-trivial on $Ker(f_{*})$ by its definition. Now Proposition 
\ref{normsubgr} implies that $\rho_{j}$ is a unitary
character. Therefore we have proved that all characters of the action of 
$\pi_{1}(C)$ on $H$ are unitary. This means that each element of the above 
action is defined by a matrix from $SL_{k}(\Z)$ with unitary eigen values. 
Applying now the theorem on the units of a ring of integers of an algebraic 
field (see, e.g. [BS, p.105,Th.2]) we obtain that each $\rho_{j}$ is a 
torsion character.\ \ \ \ \ $\Box$

Using Proposition \ref{curve} we prove the following statement.

Let $f:M\longrightarrow C$ be a holomorphic map with connected fibres of a 
compact K\"{a}hler manifold $M$ onto a smooth compact complex curve $C$ of 
genus $g\geq 1$.
\begin{Proposition}\label{pullback}
Assume that $\rho\in Hom(\pi_{1}(M),T_{2})$ and that the upper diagonal 
character
$\rho_{a}$ is pullback of a character from $Hom(\pi_{1}(C),\Co^{*})$. Assume 
also that $\rho|_{Ker(f_{*})}$ is non-trivial. Then $\rho_{a}$ is torsion.
In addition, the number of $\rho_{a}$ satisfying these assumptions is
finite.
\end{Proposition}
{\bf Proof.} First, from the assumption of the proposition it follows that 
$\rho$ maps $Ker(f_{*})$ to the abelian group $N:=DT_{2}$ of unipotent 
matrices in $T_{2}(\Co)$. Let $G$ be the intersection of kernels of all 
homomorphisms $\rho$ satisfying the conditions of the proposition. 
Then $G$ is a 
normal subgroup of $\pi_{1}(M)$ and $D(Ker(f_{*}))$ is contained in $G$. 
Set $G_{1}:=Ker(f_{*})\cap G$ and let $R:=\pi_{1}(M)/G_{1}$, 
$H:=Ker(f_{*})/G_{1}$. Then $R$ is defined by
the exact sequence 
$$
\{e\}\longrightarrow H\longrightarrow
R\longrightarrow\pi_{1}(C)\longrightarrow\{e\}\ . 
$$ 
Observe that by the definition of
$G_{1}$ and conditions of the proposition, $H$ is a non-trivial finitely 
generated free abelian group. Let us check that
$D\pi_{1}(C)$ acts trivially on $H$. In fact, let 
$a\in D\pi_{1}(M)$, $b\in Ker(f_{*})$. Then 
$\rho([a,b])=[\rho(a),\rho(b)]\subset DN=\{e\}$ for any homomorphism
$\rho$ satisfying conditions of the proposition. Thus $[a,b]\in G$. Moreover,
$[a,b]\in Ker(f_{*})$ because this is a normal subgroup. Thus
$[a,b]\in G_{1}$ and so its image in $R$ is trivial. This implies the 
required statement. Further, according to Proposition \ref{curve} all 
characters of the action $s$ of $\pi_{1}(C)$ on $H$ are torsion. Consider now 
a homomorphism $\rho$ satisfying assumptions of
the proposition. Then clearly $\rho$ determines a homomorphism
$\rho':R\longrightarrow T_{2}$ with the same image. Since according to our
assumption $\rho|_{Ker(f_{*})}$ is non-trivial, $\rho'|_{H}$ is also 
non-trivial. Thus $\rho'(H)$ is a non-trivial normal subgroup of $\rho'(R)$. 
But $\rho'(H)$ is a subgroup of unipotent matrices. Hence, the action of 
$\rho'(R)$ on $\rho'(H)$ by conjugation is defined by multiplication of 
the non-diagonal elements of $\rho'(R)$ by elements of 
$\rho_{a}(\pi_{1}(M))$. 
So $\rho_{a}(\pi_{1}(M))$ consists of
eigen values of invertible integer matrices obtained by the natural action
$s$ of $\pi_{1}(C)$ on $H$. This implies that $\rho_{a}$ is 
torsion because as we have proved, all characters of $s$ are torsion. This
argument also shows that the number of $\rho_{a}$ satisfying conditions
of the proposition is finite.
\ \ \ \  \ $\Box$
\\
{\bf 3.4.} {\bf Proof of Theorem \ref{transf}}. 
Let $V_{\rho}$ be a flat vector bundle on $M$ associated to
$\rho\in Hom(\pi_{1}(M),T_{n}(\Co))$ satisfying the 
assumptions of the theorem. Since the homomorphism 
$d(\rho)\in Hom(\pi_{1}(M),D_{n}(\Co))$ is the pullback of a homomorphism
from $Hom(\pi_{1}(C),D_{n}(\Co))$, $V_{\rho}$ is $C^{\infty}$-trivial. 
Let $V_{i}$ be a flat vector bundle associated to the 
homomorphism $\rho_{iu}''$.
Then $V_{\xi}=\oplus_{i=1}^{n}V_{i}$ is a flat vector bundle on $C$
associated to $\xi:=diag[\rho_{1u}'',...,\rho_{nu}'']\in 
Hom(\pi_{1}(C),D_{n}^{u})$. According to the results of Sections 2.1 and 2.2,
the equivalence class of $\rho$ is uniquely defined by an upper triangular
$End(f^{*}V_{\xi})$-valued $d$-harmonic 1-form $\alpha$ with a nilpotent 
(0,1)-component such that $\alpha\wedge\alpha$ represents 0 in 
$H^{2}(M,End(f^{*}V_{\xi}))$. Moreover, according to (1) $diag(\alpha_{1})$
of the (1,0)-component $\alpha_{1}$ of $\alpha$ 
is the pullback of a harmonic 1-form on $C$. Further,
$\alpha=\sum_{i,j=1}^{n}\alpha_{ij}$ where $\alpha_{ij}$ is a $d$- harmonic
1-form with values in $f^{*}(V_{i}^{*}\otimes V_{j})$.
Any pair $(f^{*}(V_{i}^{*}\otimes V_{j}),\alpha_{ij})$ determines a 
homomorphism $\rho_{ij}\in Hom(\pi_{1}(M),T_{2})$ with
upper diagonal character $f^{*}(\rho_{iu}''\otimes(\rho_{ju}'')^{-1})$.
Since by assumption (2) the latter character is not torsion, Proposition
\ref{curve} implies that $\rho_{ij}$ is the pullback of a representation
$\pi_{1}(C)\longrightarrow T_{2}$. In particular, there is a $d$-
harmonic 1-form 
$\beta_{ij}$ with values in $V_{i}^{*}\otimes V_{j}$ such that 
$f^{*}(\beta_{ij})=\alpha_{ij}$. Thus $\alpha$ is the pullback of an 
upper triangular $End(V_{\xi})$-valued harmonic 1-form
$\beta=\sum_{i,j=1}^{n}\beta_{ij}$ with a nilpotent (0,1)-component. Now 
Example \ref{ex2} 1 says that the pair $(V_{\xi},\beta)$ determines a 
representation $\widetilde\rho\in Hom(\pi_{1}(C),T_{n}(\Co))$. 
Finally, our construction and results of Section 2.2 give 
$f^{*}(\widetilde\rho)=\rho$.\ \ \ \ \ \ $\Box$\\
{\bf 3.5.} {\bf Proof of Theorem \ref{group}.}
Let $M$ be a compact K\"{a}hler manifold with fundamental group $F$.
According to the assumption of the theorem, $F$ admits 
a surjective homomorphism $p$ onto a fundamental group  $G$ of a compact 
Riemann surface of genus $g\geq 2$. (Here $p$ is the composition 
$F\stackrel{j}{\longrightarrow}G_{1}\longrightarrow G_{1}/\Z^{k}=G$.)
Then by the Siu-Beauville theorem (see, e.g.
[ABCKT, Th.2.11]) there are a holomorphic map $f:M\longrightarrow C$
with connected fibres onto a compact complex curve $C$ of genus $g'\geq g$
and a homomorphism $h:\pi_{1}(C)\longrightarrow G$ such that
$p=h\circ f_{*}$. Let $E:=f_{*}(j^{-1}(\Z^{k}))\subset\pi_{1}(C)$ be
a normal subgroup. Assume first that $\pi_{1}(C)/E$ is finite. Then 
we have $h(E)$ is a normal subgroup of a finite index in $G$. But
$h(E)=p(j^{-1}(\Z^{k}))=\{e\}$ which is wrong. Thus
$\pi_{1}(C)/E$ is infinite and so $E$ is a free group. According to the
assumption of the theorem the number of generators $r$ of $E$ is finite.
Let $S\longrightarrow C$ be a regular covering corresponding to 
$\pi_{1}(C)/E$. If $r\geq 2$ then the group $Iso(S)$ of isometries of
$S$ is finite and since $\pi_{1}(C)/E$ is infinite we have $r\leq 1$.
If $r=1$, any discrete subgroup of $Iso(S)$ is virtually cyclic and 
in particular does not act cocompactly on $S$. Thus $r=0$ which means
that $j^{-1}(\Z^{k})\subset Ker(f_{*})$, $\pi_{1}(C)\cong G$ and
$h$ is the identity homomorphism. In particular, $K\subset Ker(f_{*})$
satisfies conditions (1), (2) of Proposition \ref{curve}. This 
implies that all characters of the action of $G$ on $\Z^{k}$ are 
torsion.\ \ \ \ \ \ $\Box$
\sect{\hspace*{-1em}. Proofs of Theorems \ref{te1},
\ref{subset} and \ref{eq}.} 
Sections 4.1 and 4.2 below contain the proof of Theorem \ref{te1}.\\
{\bf 4.1.} First we prove a two-dimensional
version of the required result.

Let $T_{2}^{u}\subset T_{2} (\subset T_{2}(\Co))$ be the subgroup of matrices 
with unitary upper diagonal elements.
As before, $\rho_{a}\in Hom(\pi_{1}(M),U_{1})$ denotes the upper
diagonal character
of $\rho\in Hom(\pi_{1}(M), T_{2}^{u})$. Recall that $X_{a}$ stands for the
normalization of the image of the Albanese map
$\alpha_{X}:X\longrightarrow Alb(X)$ and 
$\alpha_{X}^{n}:X\longrightarrow X_{a}$ denotes a holomorphic map which 
covers $\alpha_{X}$.
\begin{Lm}\label{twodim}
There is a finite Galois covering $r:M_{H}\longrightarrow M$ with an
abelian Galois group $H:=H(M)$ such that for any 
$\rho\in Hom(\pi_{1}(M),T_{2}^{u})$
there is $\widetilde\rho\in Hom(\pi_{1}(M_{H a}),T_{2}^{u})$ so that
$\rho\circ r_{*}=\widetilde\rho\circ\alpha_{M_{H *}}^{n}$.
\end{Lm}
Let us introduce some definitions.
Given a character $\xi\in Hom(\pi_{1}(M),\Co^{*})$, let $\Co_{\xi}$ denote
the associated $\pi_{1}(M)$-module. We define
$\Sigma^{1}(\pi_{1}(M))$ to be the set of characters 
$\xi\in Hom(\pi_{1}(M),\Co^{*})$ such that $H^{1}(\pi_{1}(M),\Co_{\xi})$
is nonzero.\\
{\bf Theorem ([BSC]).} {\em 
There is a finite number of surjective holomorphic 
maps with connected fibres onto smooth compact complex curves 
$f_{i}:M\longrightarrow C_{i}$ and 
torsion characters $\rho_{i},\rho_{j}'\in Hom(\pi_{1}(M),\Co^{*})$ such that}
$$
\Sigma^{1}(\pi_{1}(M))=
\bigcup_{i}\rho_{i}f_{i}^{*}(Hom(\pi_{1}(C_{i}),\Co^{*}))
\cup\bigcup_{j}\{\rho_{j}'\}\ .
$$
\begin{R}\label{bsc}
{\rm
The structure of the positive dimensional components of 
$\Sigma^{1}(\pi_{1}(M))$
was discovered by Beauville [Be]. In [S1] Simpson proved that the set of zero
dimensional components of $\Sigma^{1}(\pi_{1}(M))$ consists of torsion
characters provided $M$ is a smooth projective manifold. Finally, Campana
[C1] established this fact for an arbitrary compact K\"{a}hler manifold by
reducing the question to Simpson's theorem (the projective case).}
\end{R}
{\bf Proof of Lemma \ref{twodim}.} 
We represent $Hom(\pi_{1}(M),T_{2}^{u})$ as $R_{1}\cup R_{2}
\cup R_{3}\cup R_{4}$ with $R_{i}$ to be defined below and
prove the lemma for any component of the decomposition.

It is well-known that the equivalence class of 
$\rho\in Hom(\pi_{1}(M),T_{2}^{u})$ is defined by
an element of $H^{1}(\pi_{1}(M),\Co_{\rho_{a}})$. Let $R_{1}$
consist of homomorphisms $\rho$ such that
$H^{1}(\pi_{1}(M),\Co_{\rho_{a}})=0$. Then $\rho$ is 
equivalent to a representation in the group of diagonal $2\times 2$ matrices.
Thus if $p_{t}:M_{t}\longrightarrow M$ is the Galois
covering with the Galois group $Tor(\pi_{1}(M)/D\pi_{1}(M))$ then obviously
$\rho|_{\pi_{1}(M_{t})}$ is the pullback by $\alpha_{M_{t}}$ of a 
homomorphism from $Hom(\pi_{1}(Alb(M_{t})),T_{2}^{u})$. 

Assume now that
$H^{1}(\pi_{1}(M),\Co_{\rho_{a}})\neq 0$, or equivalently,
$\rho_{a}\in\Sigma^{1}(\pi_{1}(M))$. Let $K$ be intersection
of the kernels of all characters $\rho_{i}$ and $\rho_{j}'$ from the
BSC theorem. Then $K$ is a normal subgroup of $\pi_{1}(M)$ and the
quotient group $G:=\pi_{1}(M)/K$ is finite abelian. Let
$p:M_{G}\longrightarrow M$ be the Galois covering with the Galois group
$G$. We set $\rho_{a}':=\rho_{a}|_{\pi_{1}(M_{G})}$.
Then $\rho_{a}'$ is the upper diagonal character of 
$\rho':=\rho|_{\pi_{1}(M_{G})}$. Further, according to the BSC theorem, 
there are two possibilities:\\
(1)\ $\rho_{a}'$ is trivial;\\
(2)\ there are a holomorphic surjective
map $f$ with connected fibres onto a smooth curve $C$ and 
$\rho_{1}\in Hom(\pi_{1}(C),\Co^{*})$ such that
$\rho_{a}'=\rho_{1}\circ (f\circ p)_{*}$.\\
We will denote by $R_{2}$ the set of $\rho\in Hom(\pi_{1}(M),T_{2}^{u})$ 
satisfying (1). Then the image of $\rho|_{\pi_{1}(M_{G})}$ for 
$\rho\in R_{2}$ is abelian and consists of unipotent matrices.
Therefore $\rho|_{\pi_{1}(M_{G})}$ is the pullback by $\alpha_{M_{G}}$ of a
homomorphism from $Hom(\pi_{1}(\alpha_{M_{G}}(M_{G})),T_{2}^{u})$.

Consider now case (2). Let 
$M_{G}\stackrel{g_{1}}{\longrightarrow}C_{1}\stackrel{g_{2}}{\longrightarrow}
C$ be the Stein factorization of $f\circ p$. Here $g_{1}$ is a morphism 
with connected fibres onto a smooth curve $C_{1}$ and $g_{2}$ is a finite 
morphism. Then $\rho_{a}'$ is the pullback by $g_{1}$ of 
$g_{2}^{*}(\rho_{1})\in Hom(\pi_{1}(C_{1}),\Co^{*})$.
Denote by $R_{3}(f)$ the set of $\rho\in Hom(\pi_{1}(M),T_{2}^{u})$ for which
$\rho'|_{Ker(g_{1*})}$ is non-trivial and let $R_{3}=\cup_{i}R_{3}(f_{i})$ 
where the union is taken over all maps $f_{i}$ from the BSC theorem.
According to Proposition 
\ref{pullback}, $\rho_{a}'$ for $\rho\in R_{3}(f)$ is torsion and the number 
of such characters is finite. This implies that $\rho_{a}$ is 
also a torsion character. Moreover, if $n(\rho)$ is the number of 
distinct $\rho_{a}$ for $\rho\in R_{3}(f)$ then 
$\sum_{i}\sum_{\rho\in R_{3}(f_{i})}n(\rho)<\infty$.

Finally, denote by $R_{4}(f)$ the set of homomorphisms $\rho$ for which
$\rho'|_{Ker(g_{1*})}$ is trivial and let $R_{4}=\cup_{i}R_{4}(f_{i})$.
Then for any $\rho\in R_{4}(f)$ there is 
$\xi\in Hom(\pi_{1}(C_{1}),T_{2}^{u})$ such that $\rho'=\xi\circ g_{1*}$. 

Let $F$ be the intersection of the kernels of all characters $\rho_{a}$ for 
all
$\rho\in R_{3}$ and $F_{1}:=F\cap K\cap\pi_{1}(M_{t})$. Clearly $F_{1}$ is a 
normal subgroup of $\pi_{1}(M)$ and
$H:=\pi_{1}(M)/F_{1}$ is finite abelian. Let $r:M_{H}\longrightarrow M$
be the Galois covering with Galois group $H$. Then there is a 
covering map $r_{1}:M_{H}\longrightarrow M_{G}$ such that $r=p\circ r_{1}$.
Now for $\rho\in R_{3}(f)\cup R_{4}(f)$ let
$g_{1}:M_{G}\longrightarrow C_{1}$ be as above and 
$M_{H}\stackrel{h_{1}}{\longrightarrow}C_{2}\stackrel{h_{2}}{\longrightarrow}
C_{1}$ be the Stein factorization of $g_{1}\circ r_{1}$. Here $h_{1}$ is a 
morphism with connected fibres onto a smooth curve $C_{2}$ and $h_{2}$ is a 
finite morphism. Then $\rho_{a}'':=\rho_{a}|_{\pi_{1}(M_{H})}$ is the
pullback by $h_{1}$ of $(g_{2}\circ h_{2})^{*}(\rho_{1})\in
Hom(\pi_{1}(C_{2}),\Co^{*})$ and $\rho_{a}''$ is the upper diagonal
character of $\rho'':=\rho|_{\pi_{1}(M_{H})}$. In particular, if
$\rho\in R_{4}(f)$ then $\rho''$ is pullback by $h_{1}$ of 
$\xi\circ h_{2*}\in Hom(\pi_{1}(C_{2}),T_{2}^{u})$. Since $C_{2}$ is
a smooth curve, the Albanese map 
$\alpha_{C_{2}}:C_{2}\longrightarrow Alb(C_{2})$ is an embedding. 
Below we identify $C_{2}$ with $\alpha_{C_{2}}(C_{2})$. 
Then there is a holomorphic
surjective map $a_{C_{2}}:Alb(M_{H})\longrightarrow Alb(C_{2})$ such that 
$h_{1}=a_{C_{2}}\circ\alpha_{M_{H}}$. Thus if 
$n: M_{H a}\longrightarrow
\alpha_{M_{H}}(M_{H})$ is normalization, then $\rho''$ is the pullback
by $\alpha_{M_{H}}^{n}$ of $\xi\circ (h_{2}\circ a_{C_{2}}\circ n)_{*}
\in Hom(\pi_{1}(M_{H a}),T_{2}^{u})$. \\
Assume now that $\rho\in R_{3}(f)$.
Then $\rho_{a}''$ is trivial and the image of $\rho''$ consists of
unipotent matrices (and so is abelian). This implies that $\rho''$ is 
pullback of a homomorphism from $M_{H a}$ because any such $\rho''$ is
defined by integration of a harmonic 1-form on $M_{H}$.

So we have considered all possible implications  for 
$\rho\in Hom(\pi_{1}(M),T_{2}^{u})$ and have proved
Lemma \ref{twodim} with the above defined $H$.\ \ \ \ \ $\Box$
\begin{R}\label{equiv}
{\rm 1. Let $\widetilde\rho\in Hom(\pi_{1}(M_{H a}),T_{2}^{u})$ be a 
homomorphism from Lemma \ref{twodim}, that is 
$\rho\circ r_{*}=\widetilde\rho\circ\alpha_{M_{H *}}^{n}$, and 
$\widetilde\rho_{a}$ be the upper diagonal character of $\widetilde\rho$.
Then the arguments of the lemma show that $\widetilde\rho_{a}=\exp(\eta)$ for
some $\eta\in Hom(\pi_{1}(M_{H a}),\Co)$. In particular, $\widetilde\rho_{a}$
is uniquely defined by $\rho_{a}$. Indeed, by definition
$\alpha_{M_{H *}}^{n}$ determines an isomorphism of 
$H_{1}(M_{H},\Z)/_{torsion}$ and $H_{1}(M_{H a},\Z)/_{torsion}$.
Since also $\rho_{a}\circ r_{*}=\exp(\phi)$ with 
$\phi\in Hom(\pi_{1}(M_{H}),\Co)$ (see the definition of $H$), the above
isomorphism implies that $\phi=(\eta+\phi_{1})\circ\alpha_{M_{H *}}^{n}$,
where $\phi_{1}$ is a homomorphism of $\pi_{1}(M_{H a})$ into $2\pi i\Z$. This
shows that $\widetilde\rho_{a}$ is uniquely defined.
\\
2. In fact the lemma is valid for any representation of 
$Hom(\pi_{1}(M),T_{2}(\Co))$. However, our choice of the class is 
stipulated by the following argument.

Note that any $\rho$ in Lemma \ref{twodim} is equivalently defined by
a  flat vector bundle $E$ on $M$,
$$ 
0\longrightarrow E_{1}\longrightarrow E\longrightarrow E_{2}\longrightarrow 0 
$$ 
where $E_{2}$ is a trivial flat vector bundle of complex rank-1 and 
$E_{1}$ is a flat vector bundle with unitary structure group (associated to 
the character $\rho_{a}$). It is well known that the equivalence
classes of extensions of $E_{2}$ by $E_{1}$ are in one-to-one correspondence
with the elements of the \v{C}ech cohomology group $H^{1}(M,{\bf E_{1}})$, 
where ${\bf E_{1}}$ is the sheaf of locally constant sections of $E_{1}$ 
(see, e.g., [A, Prop.2]).
In particular, by the de Rham theorem and the Hodge decomposition we obtain 
that the equivalence class of $E$ is uniquely defined by a $d$-harmonic 
1-form $\eta$ with values in $E_{1}$ (as before the Laplacian is defined by
the flat Hermitian metric on $E_{1}$). The pullback $r^{*}E$ on $M_{H}$ 
determines the homomorphism $\rho|_{\pi_{1}(M_{H})}$. Let 
$s:M_{1}\longrightarrow M_{Ha}$ be a desingularization. (Since $M_{Ha}$ is 
normal, the fibres of $s$ are connected.)
Then there is a modification $m:M_{H}'\longrightarrow
M_{H}$ and a holomorphic surjective map $\alpha_{1}:M_{H}'\longrightarrow
M_{1}$ such that the diagram}
\begin{equation}\label{factor}
\begin{array}{ccc}
M_{H}'&\stackrel{\alpha_{1}}{\longrightarrow}&M_{1}\\
m\downarrow&&\downarrow s\\
M_{H}&\stackrel{\alpha_{M_{H a}}^{n}}{\longrightarrow}&M_{H a}
\end{array}
\end{equation}
{\rm commutes. In particular,
Lemma \ref{twodim} implies that there is a flat vector bundle $F$
on $M_{1}$ lifted from $M_{H a}$ defined by
$$ 
0\longrightarrow F_{1}\longrightarrow F\longrightarrow F_{2}\longrightarrow 0 
$$ 
where $F_{2}$ is a trivial flat vector bundle of complex rank-1 and 
$F_{1}$ is a flat vector bundle with unitary structure group and 
by a $d$-harmonic 1-form $\eta_{1}$ with values in $F_{1}$ such that
$\alpha_{1}^{*}F=m^{*}(r^{*}E)$ and, hence, $\alpha_{1}^{*}(\eta_{1})=
(r\circ m)^{*}(\eta)$. Note that according to part 1 of the remark, $F_{1}$
is uniquely defined by $E_{1}$.
Clearly $\eta_{1}|_{V}=0$ and $F_{1}|_{V}$ is 
trivial for any fibre $V$ of $s$ because 
$F|_{V}$ is a trivial flat vector bundle.}
\end{R}
{\bf 4.2.} We are ready to prove Theorem \ref{te1}.
According to Proposition \ref{triang} it suffices to 
prove the theorem for representations  $\rho\in Hom(\pi_{1}(M),T_{n}(\Co))$.
Let $V_{\rho}$ be the associated flat vector bundle and
$p_{t}:M_{t}\longrightarrow M$ be a regular covering with the 
covering group $T:=Tor(\pi_{1}(M)/D\pi_{1}(M))$. The representation
$\rho|_{\pi_{1}(M_{t})}$ is defined by $p_{t}^{*}V_{\rho}$. 
Let $d(\rho)\in Hom(\pi_{1}(M),D_{n}(\Co))$ be the representation defined by 
the diagonal of $\rho$. Then clearly
$d(\rho)|_{\pi_{1}(M_{t})}\in E_{n}(M_{t})$ (see its definition in 
Introduction). Moreover, the diagonal characters $\rho_{ii}'$ of
$d(\rho)|_{\pi_{1}(M_{t})}$ satisfy $\rho_{ii}'=
\exp(\widetilde\rho_{ii}\circ p_{t *})$ for some 
$\widetilde\rho_{ii}\in Hom(\pi_{1}(M),\Co)$.
In particular, $p_{t}^{*}V_{\rho}$ is defined by
a flat connection $\omega$ on $M_{t}\times\Co^{n}$ such that 
$\omega$ is an upper triangular 1-form and $diag(\omega)$ consists of
harmonic 1-forms lifted from $M$. Further, according to the results of
Section 2.2, $p_{t}^{*}V_{\rho}$ is also defined by a flat vector bundle
$E$ on $M_{t}$ which is
a direct sum of complex rank-1 flat vector bundles 
with structure group $U_{1}$ and by an upper triangular 1-form
$\alpha+\partial h$ with values in $End(E)$ satisfying the
flatness condition. Here $\alpha$ is $d$-harmonic with a nilpotent 
(0,1)-component and $E$ is defined by the flat connection 
$diag(\omega'')-diag(\overline{\omega''})$ where $\omega''$ is the 
(0,1)-component of $\omega$. In particular, the structure of $diag(\omega)$
shows that $E=p_{t}^{*}V$ for a flat vector bundle $V$ on $M$ which
is a direct sum $\oplus_{i=1}^{n}V_{i}$ of $C^{\infty}$-trivial
complex rank-1 flat vector bundles
with structure group $U_{1}$. Further, $End(E)=\oplus_{i,j=1}^{n}E_{ij}$ with
$E_{ij}:=p_{t}^{*}V_{i}^{*}\otimes p_{t}^{*}V_{j}$.  Then 
$\alpha$ admits a decomposition $\alpha=\sum_{i,j=1}^{n}\alpha_{ij}$
where $\alpha_{ij}$ is a harmonic 1-form with values in $E_{ij}$. 
Let ${\cal H}^{1}(E_{ij})$ be the space of $E_{ij}$-valued $d$-harmonic 
1-forms on $M_{t}$. Since $E_{ij}$ is pullback of a bundle on $M$, the
finite abelian group $T$ acts on ${\cal H}^{1}(E_{ij})$. This action is
completely reducible and ${\cal H}^{1}(E_{ij})$ is isomorphic to the
direct sum $\oplus_{k=1}^{m}H_{k}$ of unitary 1-dimensional $T$-modules.
Then $\alpha_{ij}=\sum_{k=1}^{m}\alpha_{ij}^{k}$ with
$\alpha_{ij}^{k}\in H_{k}$. Let $\xi_{ij}^{k}\in 
Hom(\pi_{1}(M_{t}),T_{2}^{u})$ be a representation
constructed by the pair $(E_{ij},\alpha_{ij}^{k})$.
\begin{Lm}\label{fini}
There exists $\widetilde\xi_{ij}^{k}\in Hom(\pi_{1}(M_{H a}),T_{2}^{u})$ 
such that
$\xi_{ij}^{k}|_{\pi_{1}(M_{H})}=
(\alpha_{M_{H}}^{n})^{*}(\widetilde\xi_{ij}^{k})$.
\end{Lm}
{\bf Proof.}
According to our construction, for any $h\in T$ we have
$h^{*}(\alpha_{ij}^{k})=\rho_{k}(h)\alpha_{ij}^{k}$. Here we regard
$h$ as a biholomorphic map of $M_{t}$ and $\rho_{k}$ is the character
determining $H_{k}$. In particular, $\rho_{k}$ determines a complex rank-1 
flat vector bundle $B_{k}$ on $M$ with structure group $U_{1}$ and 
$\alpha_{ij}^{k}$ determines a $d$-harmonic 1-form 
$\widetilde\alpha_{ij}^{k}$ with values in 
$B_{k}\otimes V_{i}^{*}\otimes V_{j}$ such that 
$p_{t}^{*}(\widetilde\alpha_{ij}^{k})=\alpha_{ij}^{k}$. Let
$\xi^{k}\in Hom(\pi_{1}(M),T_{2})$ be a representation defined by
$(B_{k}\otimes V_{i}^{*}\otimes V_{j},\widetilde\alpha_{ij}^{k})$. Then 
by Lemma \ref{twodim} there is a representation 
$\widetilde\xi_{ij}^{k}\in Hom(\pi_{1}(M_{H a}),T_{2})$ such that
$\xi^{k}\circ r_{*}=\widetilde\xi_{ij}^{k}\circ\alpha_{M_{H *}}^{n}$.
But our constructions implies that 
$\xi^{k}\circ r_{*}=\xi_{ij}^{k}|_{\pi_{1}(M_{H})}$.

The lemma is proved.\ \ \ \ \ $\Box$

Further, observe that $r^{*}E$ is the pullback of a bundle $F$ 
defined on $M_{H a}$ lifted from $Alb(M)$ (because each $V_{i}$ is defined 
by a homomorphism $\exp(\eta_{i})$ with $\eta_{i}\in Hom(\pi_{1}(M),\Co)$).
Moreover, $F$ is the direct sum $\oplus_{i=1}^{n}F_{i}$ of 
$C^{\infty}$-trivial complex rank-1 flat vector bundles with 
structure group $U_{1}$. Therefore $r^{*}(End(E))$ is the pullback of 
$End(F)=\oplus_{i,j=1}^{n}F_{ij}$, $F_{ij}:=F_{i}^{*}\otimes F_{j}$ so
that $r^{*}E_{ij}=(\alpha_{M_{H}}^{n})^{*}F_{ij}$. 
Let $m:M_{H}'\longrightarrow M_{H}$ and $s,\alpha_{1}$ be the same as
in Remark \ref{equiv} (see (\ref{factor})). Then
according to Lemma \ref{fini} and Remark \ref{equiv} applied to each
$E_{ij}$ and $\alpha_{ij}^{k}$, the form  $(r\circ m)^{*}(\alpha_{ij})$ is 
the pullback by $\alpha_{1}$ of a $d$-harmonic 1-form $\beta_{ij}$ 
with values in $s^{*}F_{ij}$ (see also the uniqueness condition of part 1 of 
the remark) such that $\beta_{ij}$ vanishes on each fibre of $s$. Thus 
$(r\circ m)^{*}(\alpha)=\alpha_{1}^{*}(\beta)$, where 
$\beta=\sum_{i,j=1}^{n}\beta_{ij}$ is a 
$d$-harmonic 1-form with values in $s^{*}(End(F))$. Observe also that
the flatness condition and $\partial\overline{\partial}$-Lemma 
imply that $\alpha\wedge\alpha$ is a $d$-closed,
$d$-exact (1,1)-form. Thus from Proposition 6.1 of [Br] it follows that
$\beta\wedge\beta$ is a $d$-closed, $d$-exact matrix (1,1)-form with
values in $End(F)$. Moreover, by the definition of $\alpha$ 
the (0,1)-component of $\beta$ is nilpotent. Then from the results of 
Section 2.2 it follows that
there is a unique (up to a flat additive summand) section $h_{1}$ such that
$\beta+\partial h_{1}$ satisfies the flatness condition. Moreover, by the
same uniqueness result applied to $\alpha$ and $h$ we have
$$
\alpha_{1}^{*}(\beta+\partial h_{1})=
(r\circ m)^{*}(\alpha+\partial h)\ .
$$
In particular, the pair $(F,\beta+\partial h_{1})$ determines a
flat vector bundle $F'$ on $M_{1}$ with triangular structure group such
that $\alpha_{1}^{*}F'=(r\circ m)^{*}V_{\rho}$.
Let $V$ be a fibre of $s$ (without loss of generality we may assume that
$V$ is smooth, for otherwise we apply the arguments below to each irreducible
component of a desingularization of $V$). Then by the definition of $\beta$,
$\beta|_{V}$ is zero and $F|_{V}$ is a trivial flat vector bundle. 
From here and the flatness condition it follows that $h_{1}|_{V}$ is 
pluriharmonic and so it is constant. Thus $F'|_{V}$ is a trivial flat
vector bundle and the argument similar to that used in the proof of
Proposition \ref{normsubgr} shows that there is a  flat
vector bundle $F''$ defined on $M_{H a}$ such that $s^{*}F''=F'$.
Let $\widetilde\rho\in Hom(\pi_{M_{Ha}},T_{n}(\Co))$ be the homomorphism
associated to $F''$. Then $\widetilde\rho\circ\alpha_{M_{H *}}^{n}=
\rho\circ r_{*}$ because $m_{*}$ determines an isomorphism of fundamental
groups.

This completes the proof of Theorem \ref{te1}.\ \ \ \ \ $\Box$ 
\begin{R}\label{nilp}
{\rm Let $\rho\in Hom(\pi_{1}(M),N_{n})$ be a homomorphism into the
group of upper triangular unipotent matrices. Then the arguments of the
proof of Theorem \ref{te1} show that $Ker(\alpha_{M *}^{n})\subset Ker(\rho)$.
In fact, in this case the above flat bundle $E$ is the direct sum of trivial
complex rank-1 flat vector bundles and so all harmonic 1-forms with values
in $End(E)$ are pullbacks from $Alb(M)$. This result can be also
obtained from a result of Campana (see e.g. 
[ABCKT, Prop.3.33]).}
\end{R}
{\bf 4.3.} {\bf Proof of Theorem \ref{subset}.}
Let $\xi\in S_{n}^{\cal O}(M)$. Condition (b) in the definition of
the class $S_{n}^{\cal O}(M)$, the arguments used at the beginning
of Section 4.2 and the results of Sections 2.2 imply that 
the equivalence class of $p_{t}^{*}(\xi)$
is uniquely defined by a $d$-harmonic matrix 1-form $\alpha$ with values in a 
bundle $E$ which is a direct sum of $C^{\infty}$-trivial complex rank-1 flat 
vector bundles with structure group $U_{1}$ lifted from $M$, such that 
$\alpha\wedge\alpha$ is $d$-exact (because the flat vector bundle on $M_{t}$ 
associated to $p_{t}^{*}(\xi)$ is defined by a flat connection
with an upper triangular (0,1)-component on $M_{t}\times\Co^{n}$).
Repeating word-for-word the arguments 
used in the proof of Theorem \ref{te1} 
in the case $\xi\in Hom(\pi_{1}(M),T_{n}(\Co))$, we obtain the required
factorization of $\xi$. It remains to prove that there is a 
$\rho\in Hom(\pi_{1}(M),D_{n}^{u})$ such that
$p_{t}^{*}(\xi)\in X_{p_{t}^{*}(\rho)}$. As above (see Section 2.2), 
$p_{t}^{*}(\xi)$ is defined by a 
matrix 1-form $(\alpha'+\partial h)+\alpha''$ with values in $End(E)$
satisfying the flatness condition. Here $\alpha'$ is a holomorphic 1-form and
$\alpha''$ is an antiholomorphic nilpotent 1-form. The flatness
condition is equivalent to the equations
$$
\begin{array}{c}
(\alpha'+\partial h)\wedge(\alpha'+\partial h)=0;\\
\\
\overline{\partial}(\alpha'+\partial h)=(\alpha'+\partial h)\wedge\alpha''+
\alpha''\wedge (\alpha'+\partial h);\\
\\
\alpha''\wedge\alpha''=0\ .
\end{array}
$$
This implies that the form $z(\alpha'+\partial h)+\alpha''$ is also flat for 
any $z\in\Co$. This form determines (by iterated path integration)
a holomorphic deformation 
$\xi(z)\in Hom(\pi_{1}(M_{t}),GL_{n}(\Co))$ such that $\xi(1)=p_{t}^{*}(\xi)$
and $\xi(0)$ is defined by $End(E)$-valued antiholomorphic nilpotent 1-form 
$\alpha''$. In particular, $\xi(0)\in Hom(\pi_{1}(M_{t}),T_{n}^{u})$ 
and the diagonal $d(\xi(0))\in Hom(\pi_{1}(M_{t}),D_{n}^{u})$ of $\xi(0)$ 
determines $E$. By the definition of $E$, there is a
$\rho\in Hom(\pi_{1}(M),D_{n}^{u})$ such that $p_{t}^{*}(\rho)=d(\xi(0))$.
We complete the proof of the theorem by
\begin{Lm}\label{coinc}
Let $Y$ be a compact K\"{a}hler manifold.
For any $\eta\in Hom(\pi_{1}(Y),T_{n}(\Co))$ we have 
$X_{\eta}\subset X_{d(\eta)}$.
\end{Lm}
{\bf Proof.}
Let $X_{1}$ be an irreducible component of $X_{\eta}$ of pure dimension 
$\geq 1$. We will prove that $X_{1}\subset X_{d(\eta)}$. 
According to the results of Section 2.3 there is an open neighbourhood
$O\subset X_{1}$ of $\eta$ such that for any $w\in O$, $w\neq\eta$,
there is a holomorphic map $f:\Di\longrightarrow X_{1}$ satisfying 
$f(0)=\eta$ and $w\in f(\Di)$. Regarding $f$ as a holomorphic
deformation of $\eta$ we have 
$f(z)(g)=\sum_{i=0}^{\infty}f_{i}(g)z^{i}$, 
$g\in\pi_{1}(Y)$, such that $f(z)\in Hom(\pi_{1}(Y),GL_{n}(\Co))$,
$f_{0}=\eta$, $f_{i}$ are functions on 
$\pi_{1}(Y)$ with values in $gl_{n}(\Co)$
and the series converges uniformly on each compact of $\Di$ for
any $g\in\pi_{1}(Y)$. Let $y\in\Di$, 
$y\neq 0$, be such that $f(y)=w$. For a holomorphic diagonal matrix
function $A(z):=diag[1,z,z^{2},...,z^{n-1}]$ we define
$f_{1}(z)=A(z)^{-1}f(yz^{n-1})A(z)$. Set
$a:=(1/|y|)^{1/(n-1)}>1$. Then
a straightforward 
calculation shows that $f_{1}(z)$ is holomorphic in the disk 
$\Di_{a}:=\{z\in\Co:|z|<a\}$, $f_{1}(0)=d(\eta)$, $f_{1}(1)=w$ and for
any $z\in\Di_{a}$, $f_{1}(z)\in Hom(\pi_{1}(Y),GL_{n}(\Co))$. From here it
follows that there is an 
irreducible component $X(w)$ of $X_{d(\eta)}$ containing
$w$. This is valid for any $w\in O$.  Thus $X_{1}$ coincides with 
$X(w_{0})$ for some $w_{0}\in O$ because
$X_{1}$ is irreducible and $O\subset X_{1}$ is open.

This completes the proof of the lemma.\ \ \ \ \ $\Box$

The above lemma implies that $p_{t}^{*}(\xi)$ in Theorem \ref{subset} 
belongs to $X_{p_{t}^{*}(\rho)}$. \ \ \ \ \ $\Box$
\begin{C}\label{codef}
Let $\rho\in S_{n}(M)$ (see the definition in Introduction).
Then there is a holomorphic deformation 
$\rho(z)\in Hom(\pi_{1}(M),T_{n}(\Co))$, $z\in\Co$, of $\rho$ such that 
$\rho(1)=\rho$ and $\rho(0)\in Hom(\pi_{1}(M),D_{n}^{u})$.
\end{C}
{\bf Proof.} The proof repeats the arguments of the proof of Theorem
\ref{subset}, because in this case $\rho$ is defined by an upper
triangular flat connection $\omega$ on $M\times\Co^{n}$ and we can use the 
results of Section 2.2. Further details might be left to the reader.
\ \ \ \ \ $\Box$\\
{\bf Proof of Theorem \ref{eq}.}
Clearly we have $G_{s}^{0}(M)\subseteq G_{su}^{0}(M)$. We will prove that
$G_{su}^{0}(M)\subseteq G_{s}^{0}(M)$. 
Let $\rho\in S_{n}(M)$, that is $\rho\in Hom(\pi_{1}(M),T_{n}(\Co))$ and 
$d(\rho)\in E_{n}$. According to Corollary \ref{codef}, there is a
holomorphic deformation $\rho(z)\in Hom(\pi_{1}(M),T_{n}(\Co))$, 
$z\in\Co$, of $\rho$
such that $\rho(1)=\rho$, $\rho(0)\in Hom(\pi_{1}(M),D_{n}^{u})$ and
$\rho(z)(g)=\sum_{i=0}^{\infty}f_{i}(g)z^{i}$,
$g\in\pi_{1}(M)$ where the $f_{i}$ are upper triangular functions on
$\pi_{1}(M)$. Let $\pi_{m}:GL_{n}(\Co)[[z]]\longrightarrow GL_{n}(A_{m}[z])$
be the quotient homomorphism defined in Section 2.3. Then
$\pi_{m}(\rho(z))\in Hom(\pi_{1}(M),GL_{n}(A_{m}[z]))$ is equivalent 
to a representation $\widetilde\rho\in Hom(\pi_{1}(M),T_{p}^{u})$
(see the arguments of the proof of Proposition \ref{triang}). Therefore
$Ker(\pi_{m}(\rho(z)))\supset G_{su}^{0}(M)$ for any $m\geq 0$. Then
Proposition \ref{kern} implies that
$$
G_{su}^{0}(M)\subset\bigcap_{i\geq 0}Ker(\pi_{i}(\rho(z)))=
\bigcap_{z\in\Co}Ker(\rho(z))\subset Ker(\rho)\ .
$$
This shows that $G_{su}^{0}(M)\subset G_{s}^{0}(M)$ and completes the 
proof of the theorem.\ \ \ \ \ $\Box$
\sect{\hspace*{-1em}. Proofs of Theorem \ref{te2} and Other Results.} 
According to Theorem \ref{te1}, 
$Ker(\alpha_{M_{H *}}^{n})\subset G_{s}(M)$.
We assume that the Albanese map $\alpha_{M_{H}}$ is non-trivial. For
otherwise, $Ker(\alpha_{M_{H *}}^{n})=\pi_{1}(M_{H})$ which implies that
$M_{s}=M_{H}$ and the required statement is trivial.  Let 
$K:=Ker(\alpha_{M_{H *}})$.
Then $K$ is a normal subgroup of $\pi_{1}(M)$ (because
$\pi_{1}(M_{H})\subset\pi_{1}(M)$ is normal) and $\pi_{1}(M)/K$ is a
solvable group of the form
$$
\{e\}\longrightarrow\Z^{2l}\longrightarrow\pi_{1}(M)/K\longrightarrow
H\longrightarrow\{e\}\ .
$$
Here $l=dim_{\Co}Alb(M_{H})$. It is easy to see that $\pi_{1}(M)/K$ admits
a faithful representation to some $T_{n}(\Co)$ (see, e.g. the construction
used in Lemma \ref{fini}). Therefore
$G_{s}(M)\subset K$. Let $F:=\alpha_{M_{H *}}^{n}(G_{s}(M))$ and 
$M_{F a}\longrightarrow M_{H a}$ be a covering corresponding to the group
$F\subset\pi_{1}(M_{H a})$. Here $M_{H a}=\alpha_{M_{H}}^{n}(M_{H})$. 
By the covering homotopy theorem there is a holomorphic surjective map
$f:M_{s}\longrightarrow M_{F a}$ which covers $\alpha_{M_{H}}^{n}$.
Let $i:V\hookrightarrow M_{H}$ be a connected component of a fibre of 
$\alpha_{M_{H}}^{n}$. Then
$i_{*}(\pi_{1}(V))\subset Ker(\alpha_{M_{H *}}^{n})\subset G_{s}(M)$ and
therefore there is an embedding $i':F\hookrightarrow M_{s}$ such
that $i'(F)$ is a connected component of a fibre of $f$. This shows that $f$ 
is a proper surjective holomorphic map onto $M_{F a}$. Let 
$f=f_{1}\circ f_{2}$ be the Stein factorization of $f$.
Here $f_{2}:M_{s}\longrightarrow M'$ is a proper holomorphic map with 
connected  fibres and $f_{1}:M'\longrightarrow M_{F a}$ is a finite 
morphism. Then $M'$ is normal and any holomorphic function on $M_{s}$ is 
pullback of a holomorphic function on $M'$. We will prove that
$M'$ is holomorphically convex.

Let $K'=\alpha_{M_{H *}}^{n}(K)$ and $M_{K' a}$ be the Galois
covering of $M_{H a}$ corresponding to $K'$. Then by the definition of $K$, 
there is an intermediate covering map $p:M_{F a}\longrightarrow M_{K' a}$. 
Let $A:\Co^{l}\longrightarrow Alb(M_{H})$ be the universal covering.
Since by definition 
$\pi_{1}(M_{H})/K\cong\pi_{1}(M_{H a})/K'\cong\pi_{1}(Alb(M_{H}))$,
$M_{K' a}$ is normalization of $A^{-1}(\alpha_{M_{H}}(M_{H}))\subset\Co^{l}$.
In particular, the standard properties of Stein spaces (see e.g. [GR]) 
imply that $M_{K' a}$ is Stein and the covering $M_{F a}$ of
$M_{K' a}$ is also a normal Stein space. Moreover, $f_{1}:M'\longrightarrow
M_{F a}$ is a finite morphism and so $M'$ is a normal Stein space.
In particular, $M'$ is holomorphically convex. 

Other statements of the theorem are easy corollaries of the above argument.

The theorem is proved.\ \ \ \ \ $\Box$
\begin{R}\label{kazar}
{\rm According to Remark \ref{r-1} the kernel of any matrix unipotent 
representation of $\pi_{1}(M)$ contains $Ker(\alpha_{M *}^{n})$. Then 
the same arguments as before show that the Malcev covering of $M$ is
holomorphically convex which includes the result of Katzarkov [Ka1].}
\end{R}
{\bf Proof of Corollary \ref{stein}.} If $G_{s}(M)=\{e\}$ then 
the universal covering $M_{u}$ is holomorphically convex and any
fibre $V$ of $\alpha_{M_{H}}^{n}$ can be lifted to $M_{u}$ (we denote
the image of $V$ in $M_{u}$ by the same letter). Since $\pi_{2}(M)=0$,
Hurewicz's theorem implies that $H^{2}(M_{u},\Co)=0$. Let 
$\omega$ be a $d$-closed (1,1)-form on $M_{u}$ which is pullback of
the form on $M$ determining the K\"{a}hler metric. If $dim_{\Co}V\geq 1$ then
$\omega|_{V}$ is not $d$-exact and so $\omega$ determines a non-trivial
element of $H^{2}(M_{u},\Co)$. This contradiction shows that 
$dim_{\Co}V=0$ and the map $f_{2}:M_{u}\longrightarrow M'$ from the proof
of Theorem \ref{te2} is identity. Therefore $M_{u}=M'$ is Stein.
\ \ \ \ \ $\Box$
\begin{R}\label{cell}
{\rm
In particular, $M_{u}$ satisfying assumptions of Corollary \ref{stein}
has the structure of a $k$-dimensional complex with
$k=\dim_{\Co}M_{u}$. Thus any compact complex manifold $M$ satisfying 
the assumptions of the corollary such that $\pi_{i}(M)\neq 0$ for some
$i>dim_{\Co}M$ is not K\"{a}hler.}
\end{R}
{\bf Proof of Theorem \ref{kob}.} According to Corollary \ref{stein},
$\alpha_{M_{H}}^{n}:M_{H}\longrightarrow M_{H a}$ is a finite analytic
covering. In particular, 
$$
dim_{C}M_{H}\leq dim_{\Co}Alb(M_{H})\leq
\frac{1}{2}rank(\pi_{1}(M)/D^{2}\pi_{1}(M))\ .
$$
These inequalities and the condition of the theorem imply that
$$
dim_{C}M_{H}=dim_{\Co}Alb(M_{H})=
\frac{1}{2}rank(\pi_{1}(M)/D^{2}\pi_{1}(M))
$$
and $\alpha_{M_{H}}:M_{H}\longrightarrow Alb(M_{H})$ is a finite analytic
covering. According to Theorem \ref{subset}, any $\rho\in 
Hom(\pi_{1}(M), T_{n}(\Co))$ factors through $\alpha_{M_{H}}^{n}$
(which coincides with $\alpha_{M_{H}}$ in our case). Therefore
$\pi_{1}(M_{H})/G_{s}(M)\cong\pi_{1}(Alb(M_{H}))$.
Since $G_{s}(M)=\{e\}$, we have $\pi_{1}(M_{H})\cong\pi_{1}(Alb(M_{H}))
\cong\Z^{2l}$.
This completes the proof of the first part of the theorem.

If, in addition, $\pi_{i}(M_{u})=0$ for $1\leq i\leq dim_{\Co}M$, then
$M_{H}$ is a $K(\Z^{2l},1)$-space. Note that $Alb(M_{H})$ is also a
$K(\Z^{2l},1)$-space and $\alpha_{M_{H}}$ induces an isomorphism of 
fundamental groups. Therefore by Eilenberg's theorem (see [E]), 
$\alpha_{M_{H}}^{*}:H^{i}(Alb(M_{H}),\Z)\longrightarrow H^{i}(M_{H},\Z)$
is an isomorphism for any $i\geq 0$. Note that if 
$1\in H^{2l}(Alb(M_{H}),\Z)$ is a generator then 
$\alpha_{M_{H}}^{*}(1)$ coincides with $deg(\alpha_{M_{H}})$. Applying
the above isomorphisms with $i=2l$ we get  $deg(\alpha_{M_{H}})=1$ showing 
that $\alpha_{M_{H}}$ is biholomorphic.\ \ \ \ \ \ $\Box$

Similarly one can prove the statement of Remark \ref{rtfn}.

\end{document}